\documentclass{gtmon_a}
\pdfoutput=1
\usepackage{pb-diagram}

\proceedingstitle{Groups, homotopy and configuration spaces (Tokyo
  2005)}
\conferencestart{5 July 2005}
\conferenceend{11 July 2005}
\conferencename{Groups, homotopy and configuration spaces,
                in honour of Fred Cohen's 60th birthday}
\conferencelocation{University of Tokyo, Japan}

\editor{Norio Iwase}
\givenname{Norio}
\surname{Iwase}

\editor{Toshitake Kohno}
\givenname{Toshitake}
\surname{Kohno}

\editor{Ran Levi}
\givenname{Ran}
\surname{Levi}

\editor{Dai Tamaki}
\givenname{Dai}
\surname{Tamaki}

\editor{Jie Wu}
\givenname{Jie}
\surname{Wu}

\title[Determination of the order of the $P$--image]{Determination of the order of the $P$--image\\by Toda brackets} 

\author{Juno Mukai}                  
\givenname{Juno}
\surname{Mukai}
\address{Shinshu University\\Matsumoto 390-8621\\Japan }                  
\email{mukai@orchid.shinshu-u.ac.jp }                     
\urladdr{}

\volumenumber{13}
\issuenumber{}
\publicationyear{2008}
\papernumber{17}
\startpage{355}
\endpage{383}

\doi{}
\MR{}
\Zbl{}

\keyword{EHP sequence}
\keyword{real projective space}
\keyword{Toda bracket}
\keyword{Whitehead product}
\subject{primary}{msc2000}{55M35}
\subject{primary}{msc2000}{55Q52}
\subject{secondary}{msc2000}{57S17}

\received{24 May 2006}
\revised{21 November 2006}
\accepted{28 November 2006}
\published{19 March 2008}
\publishedonline{19 March 2008}
\proposed{}
\seconded{}
\corresponding{}
\version{}

\arxivreference{} 

\makeatletter
\def\cnewtheorem#1[#2]#3{\newtheorem{#1}{#3}[section]
\expandafter\let\csname c@#1\endcsname\c@thm}

\AtBeginDocument{\let\bar\wbar\let\tilde\wtilde\let\hat\what}

\newtheorem{mainthm}{Theorem}    
\newtheorem{thm}{Theorem}[section]    
\cnewtheorem{lem}[thm]{Lemma}          
\cnewtheorem{prop}[thm]{Proposition}
\cnewtheorem{exam}[thm]{Example}
\cnewtheorem{con}[thm]{Conjecture}

\theoremstyle{definition}
\cnewtheorem{defn}[thm]{Definition}    
\newtheorem*{rem}{Remark}             
\makeatother  
\makeautorefname{con}{Conjecture}
\newcommand{\rarrow}[1]{{\buildrel #1 \over \longrightarrow}}

\def\H{{\mathbb H}}

\def\P{{\rm P}}

\numberwithin{equation}{section}

\begin{document}

\begin{htmlabstract}
The present paper gives a proof of the author's paper in the
Proceeding of the International Conference on Homotopy Theory and
Related Topics at Korea University (2005), 109&ndash;113, on the
orders of Whitehead products of &iota;<sub>n</sub> with
&alpha;&isin;&pi;<sup>n</sup><sub>n+k</sub>, (n&ge; k+2, k&le; 24) and
improves and extends it.  The method is to use composition methods in
the homotopy groups of spheres and rotation groups.
\end{htmlabstract}

\begin{abstract} 
The present paper gives a proof of the author's paper \cite{Mu4} on
the orders of Whitehead products of $\iota_n$ with
$\alpha\in\pi^n_{n+k}$, $(n\geq k+2, k\leq 24)$ and improves and extends
it.  The method is to use composition methods in the homotopy
groups of spheres and rotation groups.
\end{abstract}

\begin{webabstract} 
The present paper gives a proof of the author's paper in the
Proceeding of the International Conference on Homotopy Theory and
Related Topics at Korea University (2005), 109--113, on the orders of
Whitehead products of $\iota_n$ with $\alpha\in\pi^n_{n+k}$, $(n\geq
k+2, k\leq 24)$ and improves and extends it.  The method is to use
composition methods in the homotopy groups of spheres and rotation
groups.
\end{webabstract}

\begin{asciiabstract}
The present paper gives a proof of the author's paper in the
Proceeding of the International Conference on Homotopy Theory and
Related Topics at Korea University (2005), 109--113, on the orders of
Whitehead products of iota_n with alpha in pi^n_{n+k}, (n > k+1, k <
25) and improves and extends it.  The method is to use composition
methods in the homotopy groups of spheres and rotation groups.
\end{asciiabstract}

\maketitle

\section*{Introduction}

This paper is a sequel to \cite{GM} by Golasi\'nski and the author in
the stable case. The methods are to use those of \cite{GM}. In
particular, the EHP sequence, the method and result of Toda
\cite[Chapter 11]{T} and the result of Nomura \cite{N} are essentially
used. Let $\pi^n_{n+k}$ denote the $2$ primary component of the
homotopy group $\pi_{n+k}(S^n)$ of the $n$ dimensional sphere
$S^n$. Let $\iota_n$ be the identity class of $S^n$ and
$\alpha\in\pi^n_{n+k}$ for $n\ge k+2$. Then our result about the order
of the Whitehead product $[\iota_n,\alpha]=P(E^{n-1}\alpha)$ is as
follows:

\begin{mainthm}[Main Theorem]\label{main} 
Let $n\ge k+2$ and $\alpha$ be an element of $\pi^n_{n+k}$. Then, the
order of the Whitehead product $[\iota_n,\alpha]$ for $n\equiv r\
(\bmod\ 8)$ with $0\leq r\leq 7$ is as given in Tables \ref{table1}
and \ref{table2} except as otherwise noted.
\end{mainthm}

\begin{table}[ht!]
\caption{}\label{table1}\vspace{5pt}

\tiny\arraycolsep 3pt\tabcolsep 3pt\centering
\begin{tabular}{||c||c|c|c|c|c|c|c|c||} \hline\hline
$\alpha\backslash r$&0&1&2&3&4&5&6&7 \\ \hline\hline
$\eta$&$2$ &$2$&$2$&$1$&$2$&$2$&$2$&$1$  \\ \hline
$\eta^2$&$2$&$2$&$1$&$1$&$2$&$2$&$1$&$1$  \\ \hline
$\nu$&$8$&$2$ &$4$&$2$&$8$&
$\begin{array}{cc}
2,&\not={\it 2^i-3}\\ 
1,&={\it 2^i-3}
\end{array}$
&$4$&$1$  \\ \hline
$\nu^2$&$2$&$2$&$2$&
$\begin{array}{cc}
2,&\not={\it 2^i-5}\\ 
1,&={\it 2^i-5}
\end{array}$
&$1$ &$1$&$2$&$1$  \\ \hline
$\sigma$&$16$&$2$&$16$&$2$&$16$ &$2$&$16$&
$\begin{array}{cc}
2,&{\it 7(16)}\\ 
1,&{\it 15(16)} 
\end{array}$
\\ \hline

$\eta\sigma$&$2$&$2$&$2$&$1$&$2$ &$2$&
$\begin{array}{cc}
1,&\not\equiv{\it 22(32)}\\
{}&\ge{\it 54}\\ 
2,&\equiv{\it 22(32)}\\
{}&\ge{\it 54} 
\end{array}$ 
&$1$\\ \hline
$\varepsilon$&$2$&$2$&$1$&$1$&$2$ &$2$&
$2$ &$1$  \\ \hline
$\bar{\nu}$&$2$&$2$&$2$&$1$&$2$ &$2$&$2$&$1$  \\ \hline
$\eta^2\sigma$&$2$&
$\begin{array}{cc}
2,&\not={\it 2^i-7}\\ 
1,&={\it 2^i-7}
\end{array}$
&$1$&$1$&$2$ &
$\begin{array}{cc}
1,&\not\equiv {\it 53(64)}\\ 
2,&\equiv {\it 53(64)}\\
{}&\ge {\it 117}
\end{array}$ 
&$1$&$1$ \\ \hline
$\eta\varepsilon$&$2$&$1$&$1$&$1$&$2$ &
$\begin{array}{cc}
1,&\not\equiv {\it 53(64)}\\
2,&\equiv {\it 53(64)}\\
{}&\ge {\it 117} 
\end{array}$ 
&$1$&$1$  
\\ \hline
$\nu^3$&$2$&
$\begin{array}{cc}
2,&\not={\it 2^i-7}\\ 
1,&={\it 2^i-7}
\end{array}$
&$1$&$1$&$1$ &$1$&$1$&$1$  \\ \hline
$\mu$&$2$&$2$&$2$&$1$&$2$ &$2$&$2$&$1$  \\ \hline
$\eta\mu$&$2$&$2$&$1$&$1$&$2$ &$2$&$1$&$1$  \\ \hline
$\zeta$&$8$&$1$&$4$&
$\begin{array}{cc}
1,&\not\equiv {\it 115(128)}\\ 
2,&\equiv {\it 115(128)}\\
{}&\ge {\it 243}
\end{array}
$
&$8$&
$1$&$4$&$1$  \\ \hline
$\sigma^2$&$2$, {\it 0(16)}&
$\begin{array}{cc}
2,&{\it 1(16)}\\ 
1,&{\it 9(16)}
\end{array}$
&$2$&
$\begin{array}{cc}
2,&{\it 3(16)}\\ 
1,&{\it 11(16)}
\end{array}$
&$2$&$2$&$2$&
$1, {\it 15(16)}$\\ \hline
$\kappa$&$2$ &$2$&$2$&$2$&$2$&$2$&$2$&$1$  \\ \hline\hline
\end{tabular}

\vspace{.3cm}
For example,\   
$\{\begin{array}{cc}
2,&\not={\it 2^i-3}\\ 
1,&={\it 2^i-3}
\end{array}\}$,
$\{ \begin{array}{cc}
2,&{\it 7(16)}\\ 
1,&{\it 15(16)} 
\end{array}\}$ \ 
and \ 
$\{ 2, {\it 0(16)}\}$ \
mean \ 
$
\{ \begin{array}{cc}
2,&\mbox{for}\ n\not={\it 2^i-3\geq 5}\\ 
1,&\mbox{for}\ n={\it 2^i-3\geq 5}
\end{array}\}$,\\ 
$\{ \begin{array}{cc}
2,&\mbox{for}\ n\equiv{\it 7\ (\bmod\ 16)\geq 23}\\ 
1,&\mbox{for}\ n\equiv{\it 15\ (\bmod\ 16)\geq 15} 
\end{array}\}$ \ and \ 
$ \{\begin{array}{cc}
2,&\mbox{for}\ n\equiv{\it 0\ (\bmod\ 16)\geq 16}\\ 
{\rm unsettled},&\mbox{for}\ n\equiv{\it 8\ (\bmod\ 16)\geq 24} 
\end{array}$\},
respectively.
\end{table}

\def\closer{\hspace{-3pt}}

\begin{table}\caption{}\label{table2}\vspace{5pt}
\tiny\arraycolsep 3pt\tabcolsep 3pt\centering
\begin{tabular}{||c||c|c|c|c|c|c|c|c||} \hline\hline
$\alpha\backslash r$&0&1&2&3&4&5&6&7 \\ \hline\hline
$\eta\kappa$&$2$&$1$&$1$&$1$&$2$ 
&$2$&$1$&$1$  \\ \hline
$\rho$&$32$&$2$ &$32$&$2$&$32$&$2$&$32$&
$a$\\ \hline
$\eta\rho$&$2$&$2$&$2$&$1$&$2$ &$2$&
{\tiny\closer
$\begin{array}{cc}
1,&\not\equiv {\it 2^9-18(2^9)}\\ 
2,&\equiv {\it 2^9-18(2^9)}\\
{} &\ge {\it 2^{10}-18}
\end{array}$\closer}
&$1$  \\ \hline
$\eta^*$&$2$&$2$&$2$&$1$&$2$ &$2$
&$2$,14(16)&$1$  \\ \hline
$\eta\eta^*$&$2$&$2$&$1$&$1$&$2$&$2$, 13(16)&$1$&$1$  \\ \hline
$\eta^2\rho$&$2$&$2$&$1$&$1$&$2$ &
{\tiny\closer $\begin{array}{cc}
1,&\not\equiv {\it 2^{10}-19(2^{10})}\\ 
2,&\equiv {\it 2^{10}-19(2^{10})}\\
{}&\ge {\it 2^{11}-19}
\end{array}$\closer}
&$1$&$1$  \\ \hline
$\nu\kappa$&$2$&$1$&$2$&$2$&$2$ &$1$&$1$&$1$  \\ \hline
$\bar{\mu}$&$2$&$2$&$2$&$1$&$2$ &$2$&$2$&$1$  \\ \hline
$\eta\bar{\mu}$&$2$&$2$&$1$&$1$&$2$ &$2$&$1$&$1$  \\ \hline
$\nu^*$&$8$&$2$&$4$&$2$&$8\ \mbox{or}\ 4$&${}$&$4$&$1$  \\ \hline
$\bar{\zeta}$&$8$&$1$&$4$&
{\tiny\closer $\begin{array}{cc}
1,&\not\equiv {\it 2^{11}-21(2^{11})}\\ 
2,&\equiv {\it 2^{11}-21(2^{11})}\\
{}&\ge {\it 2^{12}-21}
\end{array}$\closer}
&$8$&$1$
&$4$&$1$  \\ \hline
$\bar{\sigma}$&$2$&$2$&$2$&$2$&${}$&$1$, 5(16)&
$1$, 6(16)&$1$  \\ \hline
$\bar{\kappa}$&$8$&$2$&$8\ \mbox{or}\ 4$
&$2$&$4$ &$2$&$4$&$1$ \\ \hline
$\sigma^3$&$1$, 8(16)&$1$,  9(16)&$2$&$1$, 11(16)&$1$&$1$&$2$&$1$ \\ \hline
$\eta\bar{\kappa}$&$2$&$2$&$2$&$1$&$2$ &$2$&$1$&$1$  \\ \hline
$\eta^2\bar{\kappa}$&$2$&$2$&$1$&$1$&$2$ &$1$&$1$&$1$  \\ \hline
$\nu\bar{\sigma}$&$2^{\{*\}}$&${}$&${}$&$1$, 3(16)&$1$ &$1$&$1$&$1$  \\ \hline
$\eta^*\sigma$&$2$&$2$&$1$&$1$&$2$ &$2$&$1$, 6(16)&$1$  \\ \hline
$\nu\bar{\kappa}$&$8\ \mbox{or}\ 4$&${}$&$4$&$2$&$4$ &$1$&$4$&$1$ \\ \hline
$\bar{\rho}$&$16$&$2$&$16$&$2$&$16$&$2$&$16$&$b$\\ \hline
$\eta\bar{\rho}$&$2$&$2$&$2$&$1$&$2$ &$2$&
{\tiny\closer $\begin{array}{cc}
1,&\not\equiv {\it 2^{13}-26(2^{13})}\\ 
2,&\equiv {\it 2^{13}-26(2^{13})}\\
{} &\ge {\it 2^{14}-26}
\end{array}$\closer}
&$1$  \\ \hline
$\eta\eta^*\sigma$&$2$&$1$&$1$&$1$&$2$&$1$, 5(16)&$1$&$1$  \\ \hline
$\mu_{3,\ast}$&$2$&$2$&$2$&$1$&$2$ &$2$&$2$&$1$  \\ \hline
$\eta^2\bar\rho$&$2$&$2$&$1$&$1$&$2$ &
{\tiny\closer $\begin{array}{cc}
1,&\not\equiv {\it 2^{14}-27(2^{14})}\\ 
2,&\equiv {\it 2^{14}-27(2^{14})}\\
{}&\ge {\it 2^{15}-27}
\end{array}$\closer}
&$1$&$1$  \\ \hline
$\eta\mu_{3,\ast}$&$2$&$2$&$1$&$1$&$2$ &$2$&$1$&$1$  \\ \hline
$\nu^2\bar{\kappa}$&${}$&$1$&$2$&$1$
&$1$ &$1$&$2$&$1$  \\ \hline
$\zeta_{3,\ast}$&$8$&$1$&$4$&
{\tiny\closer $\begin{array}{cc}
1,&\not\equiv {\it 2^{15}-29(2^{15})}\\ 
2,&\equiv {\it 2^{15}-29(2^{15})}\\
{}&\ge {\it 2^{16}-29}
\end{array}$\closer}
&$8$&$1$
&$4$&$1$  \\ \hline\hline
\end{tabular}

\vspace{.2cm}

\hspace{-5cm}
$\{*\}$ The result holds if $\langle\bar{\nu},\sigma,\bar{\nu}\rangle
=\eta\eta^*\sigma$.

\medskip
$a=\left\{\begin{array}{ll}
1,&n\not\equiv 2^8-17(2^8);\\ 
2,&n\equiv 2^8-17(2^8)\ge 2^9-17,\ 
\end{array}\right. $
$b=\left\{\begin{array}{ll}
1,&n\not\equiv 2^{12}-25(2^{12});\\ 
2,&n\equiv 2^{12}-25(2^{12})\ge 2^{13}-25.  
\end{array}\right. 
$
\end{table}

\section[Results from Golasi\'nski and Mukai]{Results from \cite{GM}}

In this section, we shall collect the result of \cite{GM} that we
need.  We denote by $SO(n)$ the $n$-th rotation group and by $\Delta\co 
\pi_k(S^n)\to\pi_{k-1}(SO(n))$ the connecting homomorphism. The
notation $n\equiv i\ (\bmod\ k)$ is often written $n\equiv i\
(k)$. From the fact that $\pi_{4n+3}(SO(4n+3))\cong\Z$ \cite{K}, we have
$\Delta\eta_{4n+3}=0.$
\begin{gather*}
[\iota_n, \eta]=0 \ \mbox{if and only if} \ n\equiv 3\ (4)\;\mbox{or}\;n=2,6;
\tag*{\hbox{We recall}}\\
[\iota_n, \eta^2]=0 \ \mbox{if and only if} \ n\equiv 2,3\ (4)\ \mbox{or}\ n=5.
\end{gather*}
Here $\eta$ and $\eta^2$ mean exactly $\eta_n\in\pi^n_{n+1}$ and $\eta^2_n\in\pi^n_{n+2}$, respectively. Hereafter we deal with the $2$ primary components. Denote by $\sharp\alpha$ the order of $\alpha$ in a group. We recall 
$$
\sharp[\iota_n,\nu]=\left\{\begin{array}{ll}
8&\quad \mbox{if} \;n\equiv 0\ (4)\geq 8,\ n\neq 12;\\
4&\quad\mbox{if} \;n\equiv 2\ (4)\geq 6, n=4, 12;\\
2&\quad\mbox{if} \; n\equiv 1,3,5\ (8)\geq 9,\ n\neq 2^i - 3;\\
1&\quad\mbox{if}\;n\equiv 7\ (8),\ n= 2^i - 3\geq 5.
\end{array}
\right.
$$
We also recall
$$ 
\Delta(\nu^2_{8n+k})=0\ \mbox{if}\ n\geq 0\ \mbox{and}\ k=4,5.  
$$
The following is one of the main results in \cite{GM}:
\begin{thm}\label{wnu2}
$[\iota_n,\nu^2]=0 \ \mbox{if and only if} \ n\equiv 4,5,7\ (8)\ \mbox{or} \ n = 2^i - 5$ for $i\geq 4$. 
\end{thm}

Let $n\equiv 7\ (16)\geq 23$. Then, there exists an element $\delta_{n-7}\in\pi^{n-7}_{2n-8}$ satisfying
\begin{equation} \label{des7}
[\iota_n,\iota]=E^7\delta_{n-7} \ \mbox{and} \ H\delta_{n-7}=\sigma_{2n-15}\
\mbox{if} \ n\equiv 7\ (16)\geq 23.
\end{equation}
We recall
$$
\sharp[\iota_n,\sigma]=\left\{\begin{array}{ll}
16&\ \mbox{if} \ n\equiv 0\ (2)\geq 10;\\
8&\ \mbox{if} \ n=8;\\
2&\ \mbox{if}\ n\equiv 1\ (2)\geq 9,\ n\neq 11,\ n\not\equiv 15\ (16);\\
1&\ \mbox{if} \ n=11,\ n\equiv 15\ (16).
\end{array}
\right.
$$
We also recall the elements $\tau_{2n}\in\pi^{2n}_{4n}$ and 
$\bar{\tau}_{4n}\in\pi^{4n}_{8n+2}$, which are the $J$ images of the complex and symplectic characteristic elements, respectively. They satisfy the following. 

\begin{lem} \label{tau12}\

\mbox{\em(1)}\qua
$E\tau_{2n}=[\iota_{2n+1},\iota], 2\tau_{4n+2}=[\iota_{4n+2},\eta]$ and $H\tau_{2n}
=(n+1)\eta_{4n-1}${\rm;}

\mbox{\em(2)}\qua 
$E^2\bar{\tau}_{4n}=\tau_{4n+2}$ and 
$H\bar{\tau}_{4n}=\pm(n+1)\nu_{8n-1}$.
\end{lem}

\par About the group structure of the stable $k$--stem $\pi^s_k$ for
$23\leq k\leq 29$, we recall from \cite{MMO} and \cite{Od2} the following:\  $\pi^s_{23}=\{\bar{\rho},\nu\bar{\kappa},\eta^*\sigma\}\cong\Z_{16}\oplus\Z_8\oplus\Z_2$;\ 
$\pi^s_{24}=\{\eta\bar{\rho},\eta\eta^*\sigma\}\cong(\Z_2)^2$;\ $\pi^s_{25}=\{\eta^2\bar{\rho},\mu_{3,\ast}\}\cong(\Z_2)^2$;\  
$\pi^s_{26}=\{\eta\mu_{3,\ast},\nu^2\bar{\kappa}\}\cong(\Z_2)^2$;\\ 
$\pi^s_{27}=\{\zeta_{3,\ast}\}\cong\Z_8$;\ 
$\pi^s_{28}=\{\varepsilon\bar{\kappa}\}\cong\Z_2$;\ $\pi^s_{29}=0.$

By \fullref{tau12}(1) and the property of the Whitehead product, 
$$
[\iota_{4n+2},\eta\alpha]=0
\hspace{5mm} \mbox{if} \hspace{5mm} 2\alpha=0.
$$
Especially, for the elements $\beta=\nu,\zeta,\nu^*,\bar{\zeta}, \nu\bar{\kappa},\zeta_{3,\ast}$, we know the relations $4\beta=\eta^3,\eta^2\mu,\eta^2\eta^*,\eta^2\bar{\mu},\eta^3\bar{\kappa},\eta^2\mu_{3,\ast}$. By the fact that 
$H[\iota_{4n+2},2\beta]=4\beta$, we obtain 
\begin{equation}\label{4n2zeta}
\sharp[\iota_{4n+2},\beta]=4\ (\beta=\nu,\zeta,\nu^*,\bar{\zeta},\nu\bar{\kappa},\zeta_{3,\ast}).
\end{equation}
Let $n\equiv 3\ (4)\geq 7$. Then, by the fact that $\Delta\iota_n\circ\eta_{n-1}=\Delta\eta_n=0$ and $2\eta_{n-1}=0$, a Toda bracket $\{\Delta\iota_n,\eta_{n-1},2\iota\}\subset\pi_{n+1}(SO(n))$ is defined. The following result in \cite{GM} is useful to show the triviality of the Whitehead product $[\iota_n,\alpha]$:  

\begin{lem} \label{Deps0}
Let $n\equiv 3\ (4)\geq 7$. Then,

{\em (1)}\qua $\{\Delta\iota_n,\eta_{n-1},2\iota\}=0${\rm;}

{\em (2)}\qua $\Delta(E\{\eta_{n-1},2\iota_n,\alpha\})=0$, if $\alpha\in\pi_k(S^n)$ is an element satisfying $2\iota_n\circ\alpha=0$.
\end{lem}
By \fullref{Deps0},
$$
\Delta\alpha=0\ \mbox{for}\ \alpha=\varepsilon_m,\mu_m,\bar{\mu}_m,\mu_{3,m}\ (m=4n+3\geq 3);\  
\Delta\eta^*_{4n+3}=0 \ (n\geq 4)
$$
and so,
$$
[\iota_{4n+3},\alpha]=0\ \mbox{for}\ \alpha=\varepsilon,\mu,\bar{\mu},\mu_{3,\ast}\ (n\geq 0);\  
[\iota_{4n+3},\eta^*]=0\ (n\geq 4).
$$
By \cite{M2},
$$
\sharp[\iota_n,\mu]=\left\{\begin{array}{ll}
2&\ \mbox{if} \ n\equiv 0,1,2\ (4)\ge 4;\\ 
1&\ \mbox{if} \ n\equiv 3\ (4). 
\end{array}
\right.
$$
By  \cite{BJM}, \cite{FGL} and \cite{M2}, 
$$
\sharp[\iota_n,\zeta] =\left\{\begin{array}{ll}
8&\ \mbox{if} \ n\equiv 0\ (4)\ge 8;\\
4& \ \mbox{if} \ n\equiv 2\ (4)\ge 6;\\
2&\ \mbox{if} \ n\equiv 115\ (128)\ge 243;\\
1&\ \mbox{if} \ n\equiv 1\ (2)\geq 5,\ n\not\equiv 115\ (128).
\end{array}
\right.
$$ 
The results for the other elemens in the $J$--image and $\mu$--series are stated in the table. 

\section[Concerning Toda's results]{Concerning Toda's results \cite[Chapter 11]{T}}

We denote by $\P^n$ the real $n$ dimensional projective space and set $\P^n_k=\P^n/\P^{k-1}$ for $k\leq n$. Let 
$i^{m,n}_k\co  \P^m_k\hookrightarrow\P^n_k$ and $p^n_{m,k}\co  \P^n_k\to\P^n_m$ for $0\leq k\leq m\leq n$ be the canonical inclusion and collapsing maps, respectively. We set $i^n_k=i^{n-1,n}_k$ and $p^n_k=p^n_{n-1,k}$ for $k\leq n-1$. We also 
set $i^{m,n}=i^{m,n}_1$, $p^n_m=p^n_{m,1}$.
We write simply $i$ for $i^{k,n}_k, i^n_k$ and $p$ for $p^n_k$, unless otherwise stated. 

Let $i\le 4n+k-4$. We consider the exact sequence induced from a pair\break
$(E^{n-1}\P^{n+k}_n, E^{n-1}\P^{n+k-1}_n)$ \cite[(11.11)]{T}:
$$
\pi_i(E^{n-1}\P^{n+k-1}_n)\rarrow{i_*}
\pi_i(E^{n-1}\P^{n+k}_n)\rarrow{I'_k}\pi_{i+k}(S^{2n+2k-1})\rarrow{\Delta_k}\pi_{i-1}(E^{n-1}\P^{n+k-1}_n),
$$
where $I'_k$ and $\Delta_k$ are defined by the following commutative diagram:
{\small
$$
\begin{diagram}
\node{\pi_i(E^{n-1}\P^{n+k}_n)}\arrow{e,t}{j_*}
\arrow{s,r}{=}
\node{\pi_i(E^{n-1}\P^{n+k}_n,E^{n-1}\P^{n+k-1}_n)}
\arrow{e,t}{\partial}\arrow{s,r}{p_* \cong}
\node{\pi_{i-1}(E^{n-1}\P^{n+k-1}_n)}\arrow{s,r}{=}\\
\node{\pi_i(E^{n-1}\P^{n+k}_n)}\arrow{e,t}{p_*}
\arrow{s,r}{=}
\node{\pi_i(S^{2n+k-1})}\arrow{e,t}{\Delta'}
\arrow{s,r}{E^k \cong}
\node{\pi_{i-1}(E^{n-1}\P^{n+k-1}_n)}\arrow{s,r}{=}\\
\node{\pi_i(E^{n-1}\P^{n+k}_n)}\arrow{e,t}{I'_k}
\node{\pi_{i+k}(S^{2n+2k-1})}\arrow{e,t}{\Delta_k}
\node{\pi_{i-1}(E^{n-1}\P^{n+k-1}_n).}
\end{diagram}
$$}%
We denote by 
$$
\gamma_{n,k}\co  S^n\to\P^n_k
$$
the characteristic map of the $(n+1)$--cell $e^{n+1}=\P^{n+1}_k-\P^n_k$ for $k\leq n$. 
We set 
$$
\lambda_{n,k}=E^{n-1}\gamma_{n+k-1,n}.
$$
By \cite[Lemma 11.8]{T}, 
$$
\Delta_k(E^{k+1}\alpha)=\lambda_{n,k}\circ\alpha\ (\alpha\in\pi_{i-1}(S^{2n+k-2})) \hspace{5mm}  
\mbox{if} \hspace{5mm}  i\leq 4n+k-4.
$$
We denote by $\phi(s)=\sharp\{1\leq i\leq s\mid i\equiv 0,1,2,4\ (8)\}$. 
By use of \cite[Lemma 11.8, Proposition 11.9]{T}, 
we obtain:

\begin{prop} \label{nTkl}\

\mbox{\em (1)}\qua
Let $k\geq 1$ and $i\leq 4n+k-4$. Assume that 
$$
\lambda_{n,k}\circ\alpha=i_*\beta \hspace{4mm} in \hspace{4mm} \pi_{i-1}(E^{n-1}\P^{n+k-1}_n)
$$
for $\alpha\in\pi^{2n+k-2}_{i-1}$ and $\beta\in\pi^{2n-1}_{i-1}$. Then there exists an element $\delta\in\pi^{n+1}_{i+1}$ such that  $P(E^{k+3}\alpha)=E^{k-1}\delta$ and $H\delta=\pm E^2\beta$.

\mbox{\em(2)}\qua 
Let $k\geq 2, l\geq 0, n\equiv l\ (\bmod \ 2^{\phi(k)})$ and $i\leq 4n+k-4$. Assume that 
$$
\lambda_{n,k}\circ\alpha=i_*\beta \hspace{4mm} in \hspace{4mm} \pi_{i-1}(E^{n-1}\P^{n+k-1}_n)
$$ 
for $\alpha\in\pi^{2n+k-2}_{i-1}$ and $\beta\in\pi^{2n-1}_{i-1}$. Then there exists an element $\delta\in\pi^{n+1}_{i+1}$ such that  $P(E^{k+3}\alpha)=E^{k-1}\delta$ and $H\delta=\pm E^2\beta$.
\end{prop}
Although (2) is a special case of (1), it is useful in the later arguments. Hereafter \fullref{nTkl}(2) is written \fullref{nTkl}[n;k,l]. We investigate the case $4\leq k\leq 8$. 

For $n\geq 2$, we set $M^n=E^{n-2}\P^2$.
Let $\bar{\eta}_n\in[M^{n+2},S^n]\cong\Z_4$ and $\tilde{\eta}_n\in\pi_{n+2}(M^{n+1})\cong \Z_4$ for 
$n\ge 3$ be an extension and a coextension of $\eta_n$, respectively. We know the following relations in the stable groups $\{\P^2,S^0\}$ and $\pi^s_3(\P^2)$:\  
$2\bar{\eta}=\eta^2p$ and $2\tilde{\eta}=i\eta^2$.
We use the relations
$$
\bar{\eta}\tilde{\eta}=\pm 2\nu=\langle\eta,2\iota,\eta\rangle.
$$
Toda brackets are often expressed as the stable forms. 

From the fact that $E^2\P^3=M^4\vee S^5$, we take $E^2\gamma_3=2s_1\pm(E^2i^{2,3})\tilde{\eta}_3$, where $s_1\co  S^5\hookrightarrow E^2\P^3$ is the canonical inclusion. Since $E^2p^4_3\circ(E^2i^{3,4}\circ s_1)=E^4i^{1,2}$, we regard $E^2i^{3,4}\circ s_1$ as a coextension of $E^3i^{1,2}\in\pi_4(M^5)\cong\Z_2$. Set $\tilde{\imath}_5=E^2i^{3,4}\circ s_1$. Then, by the relation 
$$
2(E^2i^{3,4}\circ s_1)=\pm(E^2i^{2,4})\tilde{\eta}_3,
$$
we obtain $\pi_5(E^2\P^4)=\{\tilde{\imath}_5\}\cong\Z_8$, where $2\tilde{\imath}_5
=\pm(E^2i^{2,4})\tilde{\eta}_3$ \cite{Mu3}. We set 
$\tilde{\imath}_{n+3}=E^{n-2}\tilde{\imath}_5\in\pi_{n+3}(E^n\P^4)\cong\Z_8\ 
(n\geq 2).$ We use the relation in the stable case: 
\begin{equation}\label{2tilio}
2\tilde{\imath}=\pm i^{2,4}\tilde{\eta}.
\end{equation}
Notice that \fullref{nTkl}[n-2;2,{\it l}] for $l=2,3$ coincides with \cite[Proposition 11.10]{T} and \fullref{nTkl}[n-3;3,{\it l}] for $l=1,3$ does with \cite[Proposition 11.11]{T}, respectively. In these cases, $\lambda_{n-k,k}\in\pi_{2n-k-2}(E^{n-k-1}\P^{n-1}_{n-k})$ is taken as follows:
$$
\lambda_{n-2,2}=\left\{\begin{array}{ll}
i\eta+2\iota&(n\equiv 0\ (4));\\
i\eta&(n\equiv 1\ (4));
\end{array}\right.
$$
$$
\lambda_{n-3,3}=\left\{\begin{array}{ll}
2s_1\pm i^{2,3}\tilde{\eta}&(n\equiv 0\ (4));\\
\gamma_{5,3}\in\langle j,\eta,2\iota\rangle
&(n\equiv 2\ (4)),
\end{array}\right.
$$
where $i=E^{n-3}i^{n-1}_{n-2}$ and $j=E^{n-4}i^{n-3,n-1}_{n-3}$. By use the last part of this formula, we have $\lambda_{n-3,3}\circ\alpha=j_*\beta$ if 
$\beta\in\langle\eta,2\iota,\alpha\rangle$. So, 
\cite[Proposition 11.11.ii)]{T} is exactly interpreted as follows:   

\begin{rem}
Let $i\leq 4n-2$ and $n\equiv 3\ (4)$. Assume that $2\alpha=0$ for $\alpha\in\pi^{2n}_{i-2}$ and $\{\eta_{2n+1},2\iota,E^2\alpha\}\ni \beta$, then $P(E^7\alpha)=E^2\beta$. 
\end{rem}
 
Hereafter we use \cite[Proposition 11.11.ii)]{T} in this version. 

We use the cell structures
$$
(\mathcal{P}^4) \hspace{0.7cm}
\P^4=\P^2\cup_{\tilde{\eta}p}CM^3; \hspace{1cm} 
(\mathcal{P}^4_2) \hspace{0.7cm}
\P^4_2=S^2\cup_{\eta p}CM^3.
$$ 
By $(\mathcal{P}^4_2)$, we obtain  $\pi^s_3(\P^4_2)=\{\tilde{\imath}'\}\cong\Z_4$ and 
$\pi^s_4(\P^4_2)=\{\tilde{\imath}'\eta\}\cong\Z_2$, where $\tilde{\imath}'=p\tilde{\imath}$ and $2\tilde{\imath}'=i\eta$. Notice that   $\gamma_4=\tilde{\imath}\eta$ and  $\gamma_{4,2}=\tilde{\imath}'\eta$. 

Now, consider the case $k=4$. $\P^{n-1}_{n-4}$ has the following cell structures:
$$
\P^{n-1}_{n-4}=\left\{\begin{array}{ll}
\P^3_0=S^0\vee\P^2\vee S^3&(n\equiv 0\ (4));\\
\P^4=\P^2\cup_{\tilde{\eta}p}CM^3&(n\equiv 1\ (4));\\
\P^5_2=\P^4_2\cup_{\tilde{\imath}'\eta}e^5\ &(n\equiv 2\ (4));\\
\P^6_3=\P^5_3\cup_{\gamma_{5,3}}e^6&(n\equiv 3\ (4)).
\end{array}\right.
$$
The following cell structure is also useful:
$$
(\mathcal{P}^6_3) \hspace{1.5cm} 
\P^6_3=M^4\cup_{i\bar{\eta}}CM^5.
$$
In general, we have   
\begin{equation}\label{gam2n}
\gamma_{2n+1,k}\in\langle i,\gamma_{2n,k},2\iota\rangle.
\end{equation}
We obtain the following:
$$
\pi^s_3(\P^3_0;2)=\{\iota,\tilde{\eta},\nu\}\cong\Z\oplus\Z_4\oplus\Z_8;\ 
\pi^s_4(\P^4)=\{\tilde{\imath}\eta,i\nu\}
\cong(\Z_2)^2;
$$
$$ 
\pi^s_5(\P^5)=\{\gamma_5\}\cong\Z;\ 
\pi^s_5(\P^5_2)
=\{\gamma_{5,2},i\nu\}\cong\Z\oplus\Z_2;
$$
$$
\pi^s_5(\P^5_3)
=\{\gamma_{5,3},i'\tilde{\eta}\}\cong\Z\oplus\Z_2
$$
where 
\begin{equation}\label{gam5}
\gamma_5\in\langle i^{4,5}\tilde{\iota},\eta,2\iota\rangle,
\end{equation}
$\gamma_{5,2}\in\langle i''\tilde{\imath}',\eta,2\iota\rangle$ and $\gamma_{5,3}\in\langle i'i,\eta,2\iota\rangle\ (i'=i^5_3,i''=i^5_2)$. 
We also obtain  
$$
\pi^s_6(\P^6_3)=\{i^{4,6}_3\tilde{\eta}\eta,i\nu\}\cong(\Z_2)^2.
$$

\begin{rem}
The indeterminacy of the  bracket $\langle i'',\tilde{\imath}'\eta,2\iota\rangle$ is $\{i^{2,5}_2\nu\}+2\pi^s_5(\P^5_2)\cong\Z_2\oplus 2\Z$. Since the squaring operation $Sq^4\co  \tilde{H}^2(\P^6_2;\Z_2)\to\tilde{H}^6(\P^6_2;\Z_2)$ is trivial, we take simply $\gamma_{5,2}\in\langle i''\tilde{\imath}',\eta,2\iota\rangle$, whose indeterminacy is $2\pi^s_5(\P^5_2)$.
\end{rem}

Notice that $\P^7_4=S^4\vee M^6\vee S^7$. 
Let $s_2\co  S^7\hookrightarrow\P^7_4$ and $t\co  M^6\hookrightarrow\P^7_4$ are the canonical inclusions, 
respectively. The cell structure of $\P^n_{n-4}$ is given as follows:
$$
(\mathcal{P}^8_4) \hspace{1.5cm} 
\P^8_4=\P^7_4\cup_{\gamma_{7,4}}e^8\  
(n\equiv 0\ (8)),
$$
where 
\begin{equation}\label{gam74}
\gamma_{7,4}=2s_2\pm t\tilde{\eta}+i\nu; 
\end{equation}
$$
(\mathcal{P}^9_5) \hspace{1.5cm}
\P^9_5=\P^8_5\cup_{\gamma_{8,5}}e^9\  (\P^8_5=E^4\P^4,\ n\equiv 1\ (8)),
$$
where 
\begin{equation}\label{gam85}
\gamma_{8,5}=\tilde{\imath}\eta+i\nu; 
\end{equation} 
$$
\P^{10}_6=\P^9_6\cup_{\gamma_{5,2}+i\nu}e^{10}\ (\P^9_6=E^4\P^5_2,\ n\equiv 2\ (8));\ 
$$
$$
(\mathcal{P}^{11}_7) \hspace{2cm}
\P^{11}_7=\P^{10}_7\cup_{i\nu}e^{11}\ (\P^{10}_7=E^4\P^6_3,\ n\equiv 3\ (8));\
$$
$$
\P^4_0\ (n\equiv 4\ (8));\ 
\P^5=\P^4\cup_{\tilde{\imath}\eta}e^5\ (n\equiv 5\ (8));\
$$
$$ 
\P^6_2=\P^5_2\cup_{\gamma_{5,2}}e^6\ (n\equiv 6\ (8));\ \P^7_3=\P^6_3\vee S^7\ (n\equiv 7\ (8)).
$$  
Notice that $(\mathcal{P}^{11}_7)$  is obtained from the triviality of 
$\gamma_{10,8}\co S^{10}\to\P^{10}_8=E^8\P^2_0$. 

Let $x(n)$ be an integer such that it is odd or even according as $n$ is even or odd. Then we can set
\begin{eqnarray*}
\lambda_{n-4,4}=\left\{\begin{array}{ll}
2\iota\pm i^{n-2}_{n-3}\tilde{\eta}+x(\frac{n}{4})i\nu&(n\equiv 0\ (4));\\
\tilde{\imath}\eta+x(\frac{n-1}{4})i\nu&(n\equiv 1\ (4));\\
\gamma_{5,2}+x(\frac{n-2}{4})i\nu&(n\equiv 2\ (4));\\
x(\frac{n-3}{4})i\nu&(n\equiv 3\ (4)).
\end{array}\right.
\end{eqnarray*}

\begin{rem}
In the case $n\equiv 0\ (4)$, exactly, 
\begin{eqnarray*}
\hspace{-1cm}
\lambda_{n-4,4}=\left\{\begin{array}{ll}
2s_2\pm t\tilde{\eta}+i\nu&(n\equiv 0\ (8));\\
2s_1\pm i^{2,3}\tilde{\eta}&(n\equiv 4\ (8)).
\end{array}\right.
\end{eqnarray*}
\end{rem}

By \fullref{nTkl}, we obtain the following. 

\begin{prop} \label{nT}
Let $i\leq 4n$ and $\alpha\in\pi^{2n+2}_{i-1}$. 

\mbox{\em (1)}\qua
Let $n\equiv 0\ (\bmod \ 4)$ and assume that $\tilde{\eta}_{2n}
\circ\alpha=2\alpha=0$. Then there exists an element 
$\delta\in\pi^{n+1}_{i+1}$ such that 
$P(E^7\alpha)=E^3\delta$ and \\   $H\delta=x(\frac{n+4}{4})\nu_{2n+1}(E^2\alpha)$.
 
\mbox{\em (2)}\qua 
Let $n\equiv 1\ (\bmod \ 4)$ and assume that 
$\tilde{\imath}_{2n+1}\eta_{2n+1}\circ\alpha=0$. 
Then there exists an element $\delta\in\pi^{n+1}_{i+1}$ such that $P(E^7\alpha)=E^3\delta$ and \\ $H\delta=x(\frac{n+3}{4})\nu_{2n+1}(E^2\alpha)$.

\mbox{\em (3)}\qua
Let $n\equiv 2\ (\bmod \ 4)$ and assume that 
$E^{2n-3}\gamma_{5,2}\circ\alpha=0$. 
Then there exists an element $\delta\in\pi^{n+1}_{i+1}$ such that $P(E^7\alpha)=E^3\delta$ and \\   $H\delta=x(\frac{n+2}{4})\nu_{2n+1}(E^2\alpha)$.

\mbox{\em (4)}\qua
Let $n\equiv 3\ (\bmod \ 4)$. 
Then there exists an element $\delta\in\pi^{n+1}_{i+1}$ such that $P(E^7\alpha)=E^3\delta$ and $H\delta=x(\frac{n+3}{4})\nu_{2n+1}(E^2\alpha)$.
\end{prop}

Notice the following: 
In \fullref{nT}(1),(3), the assumptions 
$\tilde{\eta}_{2n}\alpha=0$ and   $E^{2n-3}\gamma_{5,2}\circ\alpha=0$ imply the relations  $\eta_{2n}\alpha'=0$ and $2\iota_{2n+1}\circ\alpha'=
2\alpha'=0$ respectively, where $E\alpha'=\alpha$. 

For the case $k=8$, we obtain: 

\begin{prop} \label{nT3}
Let $n\equiv l\ (\bmod \ 8)$ and $i\le 4n+4$. Let $\alpha\in\pi^{2n+6}_{i-1}$. 

\mbox{\em (1)}\qua Assume that $\pi_{2n+6}(E^{n-1}\P^{n+7}_n)
\circ\alpha=0$. Then, $P(E^{11}\alpha)$ desuspends eight dimensions. 

\mbox{\em (2)}\qua Assume that $(\pi_{2n+6}(E^{n-1}\P^{n+7}_n)-\{i\circ\sigma\})\circ\alpha=0$ for $\alpha\in\pi^{2n+6}_{i-1}$. 
Then there exists an element $\delta\in\pi^{n+1}_{i+1}$ 
such that $P(E^{11}\alpha)=E^7\delta$ and $H\delta=x\sigma_{2n+1}(E^2\alpha)$, where $x$ is even or odd according as $n\equiv l\ (\bmod \ 16)$ or $n\equiv l+8\ (\bmod \ 16)$.
\end{prop}

Hereafter \fullref{nT3}(2) is written 
\fullref{nT3}[[n;8,r]] for $r=l$ or $l+8$. We introduce some notation. If  $[\iota_n,\alpha]$ for $\alpha\in\pi^n_m$ desuspends $k$ dimensions with Hopf invariant $\theta\in\pi^{2n-2k-1}_{n+m-k-1}$, that is, if there exists an element $\delta\in\pi^{n-k}_{n+m-k-1}$ 
satisfying $E^k\delta=[\iota_n,\alpha]$ and $H\delta=\theta$, we write 
$$
H(E^{-k}[\iota_n,\alpha])=\theta.
$$
Then, immediately we obtain $P\theta=[\iota_{n-k-1},E^{-(n-k)}\theta]=0$. $\delta$ is written 
$$
\delta=\delta(\theta)=E^{-k}[\iota_n,\alpha].
$$ 
By the fact that $\sharp[\iota_n,[\iota,\iota]]=2+(-1)^n\ (n\geq 3)$ and \cite[Corollary 7.4]{BB}, $[\iota_n,\alpha\circ\beta]=[\iota_n,\alpha]\circ E^{n-1}\beta$ for $\beta\in\pi^m_l$ and so, 
\begin{equation}\label{bb}
H(E^{-k}[\iota_n,\alpha\circ\beta])=H(E^{-k}[\iota_n,\alpha])\circ E^{n-k-1}\beta.
\end{equation}
If $[\iota_n,\alpha]\neq 0$, we write
$$
H(E^{-k}[\iota_n,\alpha]_{\neq 0})=\theta.
$$
By \fullref{tau12} and by abuse of notation for $\alpha$, we obtain 
\begin{exam}\label{ex1}\

\mbox{\em (1)}\qua
$H(E^{-1}[\iota_{2n+1},\alpha])=(n+1)\eta_{4n-1}\alpha\ (\delta=\tau_{2n}\alpha), [\iota_{4n-1},\eta\alpha]=0$.

\mbox{\em (2)}\qua
$H(E^{-3}[\iota_{4n+3},\alpha])=\pm(n+1)\nu_{8n-1}\alpha\ (\delta=\bar{\tau}_{4n}\alpha),\ [\iota_{8n-1},\nu\alpha]=0$. 
\end{exam}

\noindent 
Notice that \fullref{ex1}(1) induces \cite[Proposition 11.10.ii)]{T} and \fullref{ex1}(2) does \fullref{nT}(4). 

First of all, we write up the results obtained from \cite[Proposition 11.10]{T}.

\begin{prop}\label{T10}\

\mbox{\em (1)}\qua
Let $n\equiv 0,1\ (4)$. Then, 
$H(E^{-1}[\iota_n,\alpha_1]_{\neq 0})
=\eta\alpha_1$ for $\alpha_1=\eta,\eta\sigma,\bar{\nu},\varepsilon,\mu,
\kappa,\eta\rho,\eta^*,$ $\bar{\mu},\eta\bar{\kappa},\eta^*\sigma,\mu_{3,\ast}$ and $H(E^{-1}[\iota_n,\alpha_2])=0$ for $\alpha_2=\eta\varepsilon,\eta^2\sigma,\sigma^2,\eta\kappa,\eta^2\rho,\bar{\sigma},\nu\bar{\sigma},\eta\eta^*\sigma,$
$\eta^2\bar{\rho}$.

\mbox{\em (2)}\qua
$H(E^{-1}[\iota_{4n},\beta]_{\neq 0})=\eta\beta$ for $\beta=\eta^2,\eta\mu,\eta\eta^*,\eta\bar{\mu}, \eta^2\bar{\kappa},\eta\mu_{3,\ast}$.

\mbox{\em (3)}\qua
$H(E^{-1}[\iota_{4n+1},\delta_1]_{\neq 0})=\eta\delta_1$ for $\delta_1=\sigma,\rho,\bar{\kappa},\bar{\rho}$ and 
$H(E^{-1}[\iota_{4n+1},\delta_2])=0$ for $\delta_2=\nu,\zeta,
\nu^*,\bar{\zeta},\nu\bar{\kappa},\zeta_{3,\ast}$.

\mbox{\em (4)}\qua
If $\sharp[\iota_{4n},\nu^*]=8$, then $[\iota_{4n+1},\eta\eta^*]\neq 0$. 
\end{prop}
\begin{proof}
We prove (1) for $\kappa$. By \cite[Proposition 11.10]{T}, $H(E^{-1}[\iota_n,\kappa])=\eta\kappa$. 
Assume that $[\iota_n,\kappa]=0$. Then, by the EHP sequence, $\delta\in P\pi^{2n-1}_{2n+14}=\{[\iota_{n-1},\eta\kappa],[\iota_{n-1},\rho]\}$ for $\delta=E^{-1}[\iota_n,\kappa]$. Applying the Hopf homomorphism $H\co \pi^{n-1}_{2n+12}\to\pi^{2n-3}_{2n+12}$ to this relation implies 
$\eta\kappa=0$ for $n\equiv 0\ (\bmod\ 4)$ and 
$\eta\kappa\in\{2\rho\}$ for $n\equiv 1\ (\bmod\ 4)$. This is a contradiction. 

Next, we prove (2) for $\eta\eta^*$. Let $n\equiv 0\ (4)$.  By \cite[Proposition 11.11]{T}, $H(E^{-1}[\iota_n,\eta\eta^*])$ $=\eta^2\eta^*=4\nu^*$. The assumption $[\iota_n,\eta\eta^*]=0$ induces $\delta\in P\pi^{2n-1}_{2n-19}$ and a contradictory relation $4\nu^*=0$ for $\delta=E^{-1}[\iota_n,\eta\eta^*]$. The proof of (3) is similarly obtained. 

Finally, we show (4). Assume that $[\iota_{4n+1},\eta\eta^*]=0$. 
From the fact that  $[\iota_{4n+1},\eta\eta^*]=E(\tau_{4n}\eta\eta^*)$ and the assumption $\sharp[\iota_{4n},\nu^*]=8$, we have $\tau_{4n}\eta\eta^*\in\{4[\iota_{4n},\nu^*],[\iota_{4n},\eta\bar{\mu}]\}$. This implies a contradictory relation $4\nu^*=0$, and hence (4) follows. 
 \end{proof}

\medskip

Hereafter, ``the assumption $[\iota_n,\alpha]=0$'' is written ``${\mathcal ASM}[\alpha]$'' and ``a contradictory relation $\beta\in B$'' is written ``${\mathcal CDR}[\beta\in B]$''. 
As an application of \cite[Proposition 11.11]{T}, we show: 

\begin{prop}\label{T11}\

\mbox{\em (1)}\qua
$H(E^{-2}[\iota_{4n+2},\alpha])\in\langle\eta,2\iota,\alpha\rangle$ if $2\alpha=0$,\\ 
$H(E^{-2}[\iota_{4n+2},\alpha_1]_{\neq 0})\in\langle\eta,2\iota,\alpha_1\rangle$ for  $\alpha_1=\nu^2,8\sigma,
\sigma^2,16\rho,\sigma^3,
8\bar{\rho},\nu^2\bar{\kappa}$ and\\ 
$H(E^{-2}[\iota_{4n+2},\alpha_2])=0$ for  $\alpha_2=\eta\sigma,\bar{\nu},\varepsilon,\nu^3,
\eta\rho,\bar{\sigma},\eta\bar{\rho}$.

\mbox{\em (2)}\qua
$H(E^{-2}[\iota_{4n},\beta_1]_{\neq 0})\in\langle 2\iota,\eta,\beta_1\rangle$ for $\beta_1=\eta\kappa,\eta^2\rho, \eta\eta^*\sigma$. 

\mbox{\em (3)}\qua
$H(E^{-2}[\iota_{4n},\beta_2])=0$ for $\beta_2=4\nu,8\sigma,4\zeta,\sigma^2,16\rho,4\bar{\zeta},\bar{\sigma},4\bar{\kappa},4\nu\bar{\kappa},8\bar{\rho},4\zeta_{3,\ast}$.  
\end{prop}
\begin{proof}
Let $n\equiv 2\ (4)$. The first part of (1) is a direct consequence of \cite[Proposition 11.11.ii)]{T}. By the fact that $\langle\eta,2\iota,\sigma^2\rangle\ni\eta^*\ (\bmod\ \eta\rho)$ and \cite[Proposition 11.11.ii)]{T},
$$
H(E^{-2}[\iota_n,\sigma^2])=\eta^*.
$$
${\mathcal ASM}[\sigma^2]$ induces $E\delta\in P\pi^{2n-1}_{2n+14}
=\{[\iota_{n-1},\alpha]\}=\{E(\tau_{n-2}\alpha)\}
$ (\fullref{tau12}(1)) and
$\delta\ (\bmod\ \tau_{n-2}\rho,\tau_{n-2}\eta\kappa)
\in P\pi^{2n-3}_{2n+13}$, where $\delta
=E^{-2}[\iota_n,\sigma^2]$ and $\alpha=\rho, \eta\kappa$. Hence, ${\mathcal CDR}[\eta^*\ (\bmod\ \eta\rho)=0]$ and the second part of (1) for $\sigma^2$ follows. Next we prove the second part of (1) for $\nu^2\bar{\kappa}$. By the fact that 
$\langle\eta,2\iota,\nu^2\rangle\ni\varepsilon\ (\bmod\ \eta\sigma)$ and \cite[Proposition 11.11.ii)]{T},
$H(E^{-2}[\iota_n,\nu^2])=\varepsilon$ and $
H(E^{-2}[\iota_n,\nu^2\bar{\kappa}])
=\varepsilon\bar{\kappa}$ by \eqref{bb}. 
${\mathcal ASM}[\nu^2\bar{\kappa}]$ induces  $E(\delta\bar{\kappa})\in\{[\iota_{n-1},\zeta_{3,\ast}]\}$ and $\delta\bar{\kappa}\ (\bmod\ \tau_{n-2}\zeta_{3,\ast})\in P\pi^{2n-3}_{2n+25}$,  where $\delta=E^{-2}[\iota_n,\nu^2]$. By the relation $\eta\zeta_{3,\ast}=0$, we obtain\\ 
${\mathcal CDR}[\varepsilon\bar{\kappa}=0]$. 

The third part of (1) follows from \cite[Proposition 11.11.ii)]{T} and the fact that $\langle\eta,2\iota,\alpha_2\rangle\ni0$. By \cite[Proposition 11.11.i)]{T}, 
\begin{equation}\tag{$\diamond$} H(E^{-2}[\iota_{4n},\eta\kappa])=\langle 2\iota,\eta,\eta\kappa\rangle=\nu\kappa.\label{diamond}
\end{equation}
${\mathcal ASM}[\eta\kappa]$ implies  $E\delta\in P\pi^{2n-1}_{2n+15}
=\{E(\tau_{n-2}\eta\rho),E(\tau_{n-2}\eta^*)\}$ and\\ 
$\delta\ (\bmod\ \tau_{n-2}\eta\rho,\tau_{n-2}\eta^*)
\in P\pi^{2n-3}_{2n+14}$, where $\delta
=E^{-2}[\iota_n,\eta\kappa]$.
Hence, \\
${\mathcal CDR}[\nu\kappa\ (\bmod\ \eta^2\rho,\eta\eta^*)=0]$ and the first part of (2) follows.
By the parallel argument, the rest of the assertion follows. We use the following facts:
$\langle 2\iota,\eta,\beta_2\rangle=0$;\ 
$\langle\eta,2\iota,16\rho\rangle\ni\bar{\mu}\ (\bmod\ 
\eta^2\rho,\eta\eta^*);\ \langle\eta,2\iota,\sigma^3\rangle\ni\eta^*\sigma\ 
(\bmod\ \eta^3\bar{\kappa});\
\langle 2\iota,\eta,\eta^2\rho\rangle\ni\bar{\zeta}\ (\bmod\ 2\bar{\zeta});\ \langle 2\iota,\eta,\eta\eta^*\sigma\rangle=\nu^2\bar{\kappa}$\ \cite{HiM}.
\end{proof}

By \fullref{T11}(2), we obtain 
\begin{equation}\label{nukap5}
[\iota_{4n+1},\nu\kappa]=0
\end{equation}
$$
[\iota_{4n+1},\nu^2\bar{\kappa}]=0.\leqno{\hbox{and}} 
$$
Here we summarize Toda brackets in $\pi^s_*(\P^2)$ needed in the subsequent arguments. 
Since $\pi^s_7(\P^2)=\{i\nu^2\}\cong\Z_2$ and $\pi^s_5(\P^2)=\{\tilde{\eta}\eta^2\}\cong\Z_2$, the indeterminacy of the bracket $\langle i\bar{\eta},\tilde{\eta},\nu\rangle\subset\pi^s_8(\P^2)$ is $i\bar{\eta}\circ\pi^s_7(\P^2)+\pi^s_5(\P^2)\circ\nu=0$.
We set $\widetilde{\nu^2}=\langle i\bar{\eta},\tilde{\eta},\nu\rangle$, which is a coextension of $\nu^2$. Let $\widetilde{\sigma^2}\in\langle i,2\iota,\sigma^2\rangle\subset\pi^s_{16}(\P^2)$ be a 
coextension of $\sigma^2$ and $\overline{i\nu}\in\{M^5,\P^2\}$ an extension of $i\nu\in\pi^s_4(\P^2)$.  Then, we show:
\begin{lem} \label{coexnu2}\

\mbox{\em (1)}\qua
$\langle i\bar{\eta},\tilde{\eta},
\nu^*\rangle\ni\widetilde{\sigma^2}\sigma\ (\bmod\ i\eta^2\bar{\kappa},i\nu\bar{\sigma})$.

\mbox{\em (2)}\qua
$\langle i\nu,2\iota,\sigma^2\rangle=i\nu^*$.

\mbox{\em (3)}\qua
$\langle i\nu,2\iota,16\rho\rangle=i\bar{\zeta}$. 

\mbox{\em (4)}\qua
$\langle i\nu,2\iota,\eta^*\rangle=0$. 

\mbox{\em (5)}\qua
$\langle \overline{i\nu},\tilde{\eta},4\iota\rangle=\pi^s_7(\P^2)$. 

\mbox{\em (6)}\qua
$\langle\tilde{\eta}p,\tilde{\eta}\eta^2,\eta\rangle=0$.

\mbox{\em (7)}\qua
$\langle\tilde{\eta}p,\tilde{\eta}\eta^2,\sigma^2\rangle\ni 0\ (\bmod\ \tilde{\eta}\eta\bar{\mu})$.

\mbox{\em (8)}\qua
$\langle i\eta\bar{\eta},\tilde{\eta},\nu\rangle=
\widetilde{\nu^2}\eta=i\varepsilon$, 
$\widetilde{\nu^2}\sigma=0$ \ and \  $\widetilde{\nu^2}=\langle\tilde{\eta},\nu,\eta\rangle$.

\mbox{\em (9)}\qua
$\langle\tilde{\eta},\nu,\nu^3\rangle=i\eta\kappa$. 

\mbox{\em (10)}\qua
$\widetilde{\nu^2}\eta\eta^*=i\eta\eta^*\sigma$ \ and \ 
$\langle\tilde{\eta}p,\tilde{\eta}\eta^2,\nu^*\rangle\ni i\eta\eta^*\sigma\ (\bmod\ \tilde{\eta}\eta^2\bar{\kappa})$. 
\end{lem}
\begin{proof}
Since $\langle p,i\bar{\eta},\tilde{\eta}\rangle=\pm\nu$ and $\nu\nu^*=\sigma^3$, 
we have $p\circ\langle i\bar{\eta},\tilde{\eta},
\nu^*\rangle=\sigma^3$. This leads to (1). By the fact that $\nu^*\in\langle\nu,2\sigma,\sigma\rangle$ and $\nu\circ\pi^s_{15}=0$, we see that 
$$
\langle i\nu,2\iota,\sigma^2\rangle\subset\langle i\nu,2\sigma,\sigma\rangle\ni i\nu^*\ (\bmod\ i\nu\circ\pi^s_{15}+\pi^s_{12}(\P^2)\circ\sigma=\{\tilde{\eta}\mu\sigma\}).
$$
We have $p\circ\langle i\nu,2\iota,\sigma^2\rangle=\langle p,i\nu,2\iota\rangle\circ\sigma^2\subset\pi^s_3\circ\sigma^2=0$, $p(i\nu^*)=0$ and $p(\tilde{\eta}\mu\sigma)=\eta\mu\sigma=\eta^2\rho$. This leads to (2).

We obtain 
$$
\langle i\nu,2\iota,16\rho\rangle\subset\langle i\nu,8\iota,4\rho\rangle\supset i\circ\langle \nu,8\iota,4\rho\rangle\ni i\bar{\zeta}
$$
$$
(\bmod\ 
i\nu\circ\pi^s_{16}+\pi^s_5(\P^2)\circ 4\rho=0).
$$
We get that  
$$
\langle i\nu,2\iota,\eta^*\rangle\subset\langle i,2\nu,\eta^*\rangle\supset\langle i,2\iota,0\rangle\ni 0\  (\bmod\ i_*\pi^s_{20}+\pi^s_5(\P^2)\circ\eta^*).
$$
Since $\tilde{\eta}\eta^2\eta^*=4\tilde{\eta}\nu^*=0$, the indeterminacy is $\{i\bar{\kappa}\}$. Hence, (4) follows from the fact that 
$\langle\bar{\eta},i\nu,2\iota\rangle\subset\pi^s_5=0$ and $\bar{\eta}\circ i\bar{\kappa}=\eta\bar{\kappa}$.

The indeterminacy of 
$\langle\overline{i\nu},\tilde{\eta},4\iota\rangle$ contains $\overline{i\nu}\circ\pi^s_4(\P^2)=\{i\nu^2\}=
\pi^s_7(\P^2)$. 

We obtain 
$$
\langle\tilde{\eta}p,\tilde{\eta}\eta^2,\eta\rangle\subset
\langle\tilde{\eta},4\nu,\eta\rangle\supset\langle 0,\nu,\eta\rangle\ni 0\ (\bmod\ \tilde{\eta}\circ\pi^s_5
+\pi^s_7(\P^2)\circ\eta=0).
$$
We see that 
$$
\langle\tilde{\eta}p,\tilde{\eta}\eta^2,\sigma^2\rangle\subset
\langle\tilde{\eta},4\nu,\sigma^2\rangle\ni 0\ (\bmod\ \tilde{\eta}\circ\pi^s_{18}+\pi^s_7(\P^2)\circ\sigma^2),
$$
where $\pi^s_7(\P^2)\circ\sigma^2=0$ and $\tilde{\eta}\nu^*=0$ because $\langle 2\iota,\eta,\nu^*\rangle\subset\{2\bar{\kappa}\}$. This leads to (7). 

By the equality 
$\langle\eta\bar{\eta},\tilde{\eta},\nu\rangle=\varepsilon$ \cite[Lemma 4.2]{GM}, $i\bar{\eta}\widetilde{\nu^2}\in i\langle\eta\bar{\eta},\tilde{\eta},\nu\rangle=i\varepsilon.
$ This implies $\bar{\eta}\widetilde{\nu^2}=\varepsilon$.
We have $i\varepsilon\in\langle i\eta\bar{\eta},\tilde{\eta},\nu\rangle\ (\bmod\ i\eta\bar{\eta}\circ\pi^s_7(\P^2)+\pi^s_6(\P^2)\circ\nu=0)$ and $\widetilde{\nu^2}\eta\in i\langle 2\iota,\nu^2,\eta\rangle\ni i\varepsilon\ (\bmod\ i\eta\sigma).$
Composing $\bar{\eta}$ on the left to this relation yields  $\widetilde{\nu^2}\eta=i\varepsilon$. 
We have $\widetilde{\nu^2}\sigma=
\langle i\bar{\eta},\tilde{\eta},\nu\rangle\circ\sigma= i\bar{\eta}\circ\langle\tilde{\eta},\nu,\sigma\rangle=0.$
Since $p\circ\langle\tilde{\eta},\nu,\eta\rangle=\nu^2$, we can set $\langle\tilde{\eta},\nu,\eta\rangle=\widetilde{\nu^2}+ai\sigma$ for $a\in\{0,1\}$. By the fact that 
$\eta\bar{\eta}\circ\langle\tilde{\eta},\nu,\eta\rangle
=\langle\eta\bar{\eta},\tilde{\eta},\nu\rangle\circ\eta
=\eta\varepsilon$ and $\eta\bar{\eta}(\widetilde{\nu^2}+ai\sigma)=\eta\varepsilon+a\eta^2\sigma$, we have $a=0$. 

By the relations $\nu^3=\eta\bar{\nu}$, $\langle 2\iota,\nu^2,\bar{\nu}\rangle\ni\eta\kappa\ (\bmod\ 2\rho)$ and (8),  
$$
\langle\tilde{\eta},\nu,\nu^3\rangle\supset\langle\tilde{\eta},\nu,\eta\rangle\circ\bar{\nu}=\widetilde{\nu^2}\bar{\nu}\in i\langle 2\iota,\nu^2,\bar{\nu}\rangle=i\eta\kappa
$$
$$ 
(\bmod\ \tilde{\eta}\circ\pi^s_{13}
+\pi^s_7(\P^2)\circ\nu^3=0).
$$
By (8) and \cite[(6.3)]{MMO}, $\widetilde{\nu^2}\eta\eta^*
=i\varepsilon\eta^*=i\eta\eta^*\sigma.$
By the fact that  $2\tilde{\eta}\bar{\eta}=\tilde{\eta}\eta^2p=i\eta\bar{\eta}\circ i\bar{\eta}$, $p_*\pi^s_{24}(\P^2)=\pi^s_{22}=\{\eta^2\bar{\kappa},\nu\bar{\sigma}\}\cong(\Z_2)^2$ and (1),   
$$
\langle\tilde{\eta}p,\tilde{\eta}\eta^2,\nu^*\rangle
\supset\langle\tilde{\eta}\eta^2p,\tilde{\eta},\nu^*\rangle\supset i\eta\bar{\eta}\circ\langle
i\bar{\eta},\tilde{\eta},\nu^*\rangle
\ni i\eta\eta^*\sigma
$$
$$
(\bmod\ \tilde{\eta}p\circ\pi^s_{24}(\P^2)
+\pi^s_7(\P^2)\circ\nu^*=\{\tilde{\eta}\eta^2\bar{\kappa}\}).
$$
This leads to (10). 
\end{proof}

We recall from \cite{Mu2} that $\{\P^4,S^0\}=\{\bar{\eta}'\}\cong\Z_8$ and
\begin{equation}\label{etatild}
\bar{\eta}'\tilde{\imath}=\nu, \hspace{5mm} \mbox{where} \hspace{5mm} \bar{\eta}'\in\langle\bar{\eta},\tilde{\eta}p,p_{4,2}\rangle. 
\end{equation}
We obtain the following. 

\begin{lem}\label{P76}\

\mbox{\em(1)}\qua
$\pi^s_7(\P^4)=\{\widetilde{\tilde{\eta}\eta^2},i\nu^2\}\cong(\Z_2)^2$ and $\pi^s_7(\P^6)=\{\tilde{\eta}',i\nu^2\}
\cong\Z_8\oplus\Z_2$, where $\widetilde{\tilde{\eta}\eta^2}\in\langle i^{2,4},\tilde{\eta}p,\tilde{\eta}\eta^2\rangle$,  $\tilde{\eta}'\in\langle i^{4,6},\tilde{\imath}\bar{\eta},\tilde{\eta}\rangle$ and $4\tilde{\eta}'\equiv i^{4,6}\widetilde{\tilde{\eta}\eta^2}\ (\bmod\ i\nu^2)$. 

\mbox{\em(2)}\qua
 $\pi^s_7(\P^6_3)=\{\tilde{\eta}''\}\cong\Z_8$, where $\tilde{\eta}''=p^6_3\tilde{\eta}'$.

\mbox{\em(3)}\qua
$\pi^s_7(\P^4)\circ\eta=\pi^s_7(\P^4)\circ\sigma^2=0$\  and \ $\pi^s_7(\P^6)\circ\sigma^2
=\pi^s_7(\P^6_3)\circ\sigma^2=0$.
\end{lem}
\begin{proof}
(1) is just \cite[Proposition 4.1]{Mu2}. (2) is obtained by use of the cell structure $(\mathcal{P}^6_3)$ and (1). The first two equalities in (3) are obtained by \fullref{coexnu2}(6),(7) and the relation $i^{2,4}\tilde{\eta}p=0\in\{M^3,\P^4\}$. 
To show the next two equalities in (3), it suffices to prove $\langle\tilde{\imath}\bar{\eta},\tilde{\eta},\sigma^2\rangle\ni 0$.  By \eqref{2tilio}, the relation  $\langle\tilde{\eta},\nu,\sigma\rangle=0$ and the second equality in (3),  
$$
\langle\tilde{\imath}\bar{\eta},\tilde{\eta},\sigma^2\rangle\subset\langle\tilde{\imath},2\nu,\sigma^2\rangle
\supset\langle i^{2,4}\tilde{\eta},\nu,\sigma\rangle\circ\sigma
\ni 0\ (\bmod\ \tilde{\imath}\circ\pi^s_{18}).
$$ 
We have $2\tilde{\imath}\nu^*=i^{2,4}\tilde{\eta}\nu^*=0$. 
By the fact that $\{M^6,S^0\}=\{\nu^2p\}\cong\Z_2$, \eqref{etatild} and (1), $\bar{\eta}'\circ\tilde{\imath}\nu^*=\sigma^3$, $\bar{\eta}'\circ\langle\tilde{\imath}\bar{\eta},\tilde{\eta},\sigma^2\rangle=\langle\bar{\eta}',\tilde{\imath}\bar{\eta},\tilde{\eta}\rangle\circ\sigma^2$ and   $8\langle\bar{\eta}',\tilde{\imath}\bar{\eta},\tilde{\eta}\rangle=\langle 8\iota,\bar{\eta}',\tilde{\imath}\bar{\eta}\rangle\circ\tilde{\eta}\subset\{M^6,S^0\}\circ\tilde{\eta}=0$. This implies  $\langle\bar{\eta}',\tilde{\imath}\bar{\eta},\tilde{\eta}\rangle\subset 2\pi^s_7$ and $\bar{\eta}'\circ\langle\tilde{\imath}\bar{\eta},
\tilde{\eta},\sigma^2\rangle=0$.
\end{proof}

We show:
\begin{lem}\label{nTa1}\

\mbox{\em (1)}\qua
$H(E^{-3}[\iota_{4n},\alpha])=\frac{1+(-1)^n}{2}\nu\alpha$ for $\alpha=4\nu,\nu^2,8\sigma,
\nu^3,4\zeta,16\rho$,\\
$\nu\kappa,4\nu^*,4\bar{\zeta},\bar{\sigma},4\bar{\kappa},4\nu\bar{\kappa},8\bar{\rho}$. In particular,
$H(E^{-3}[\iota_{8n},\alpha])=\nu\alpha$\ for $\alpha=\nu^2,\nu\kappa$,\\
$\bar{\sigma},4\bar{\kappa}$.

\mbox{\em (2)}\qua
$H(E^{-7}[\iota_{8n},\beta])=0$ for  $\beta=8\sigma,16\rho,8\bar{\rho}$.

\mbox{\em (3)}\qua
$H(E^{-7}[\iota_{8n},\sigma^2])=0$ or $\sigma^3$.
\end{lem}
\begin{proof}
(1) is a direct consequence of \fullref{nT}(1). Let $n\equiv 0\ (8)$. We have $\P^{n-1}_{n-8}
=E^{n-8}\P^7_0$ and $\gamma_{n-1,n-8}\in 2\pi^s_7(S^7)\oplus\pi^s_7(\P^6)\oplus\pi^s_7$. 
By \fullref{P76}(1), $\lambda_{n-8,8}\circ\beta=0$. Hence, by \fullref{nT3}[[n-8;8,0]], $[\iota_n,\beta]$ desuspends eight dimensions. Similarly, by \fullref{P76}(3) and \fullref{nT3}[[n-8;8,0]], $\lambda_{n-8,8}\circ\sigma^2=0$ or $i\sigma^3$. 
\end{proof}
 
By \fullref{nTa1}(1), 
 \begin{equation}\label{nbs4}
 [\iota_{8n+4},\nu\bar{\sigma}]=0.
 \end{equation} 
 We need the following. 
\begin{lem}\label{6bamu}
$H(E^{-5}[\iota_{8n+6},\alpha])=0$ for $\alpha=\eta,
\varepsilon,\bar{\nu},\mu,\kappa,\eta^*,\nu\kappa,\bar{\mu},\bar{\sigma},\eta\bar{\kappa},\nu\bar{\sigma}$,\\ 
$\mu_{3,\ast}$.
\end{lem}
\begin{proof}
We show the assertion for $\alpha=\eta,
\varepsilon,\mu,\kappa,\eta^*,\bar{\sigma}$. Let $n\equiv 6\ (8)$. In \fullref{nTkl}[n-6;6,0], 
$\P^{n-1}_{n-6}=E^{n-6}\P^5_0$. 
We take $\lambda_{n-6,6}=\gamma_5$. 
By \eqref{gam5} and \eqref{2tilio}, $\gamma_5\eta=\pm i^{2,5}\tilde{\eta}\nu=0$. 
We obtain $\gamma_5\varepsilon=0$, because $\langle\eta,2\iota,\varepsilon\rangle=\{\eta\varepsilon\}$. By the fact that $\langle\eta,2\iota,\mu\rangle=\pm 2\zeta$ and $\langle2\iota,\eta,\zeta\rangle=0$, $$\gamma_5\mu\in i^{4,5}\tilde{\imath}\circ\langle\eta,2\iota,
\mu\rangle=i^{2,5}\tilde{\eta}\zeta=0.$$
By the relation $\langle\eta,2\iota,\eta^*\rangle\ni\pm 2\nu^*\ (\bmod\ \eta\bar{\mu})$, we have $\gamma_5\eta^*=i^{2,5}\tilde{\eta}\nu^*=0$. By the fact that $\langle\eta,
2\iota,\kappa\rangle\ni 0\ (\bmod\ \eta\rho)$ and 
$\langle\eta,
2\iota,\bar{\sigma}\rangle\ni 0\ (\bmod\ \eta\bar{\kappa})$, 
we obtain $\gamma_5\kappa
=\gamma_5\bar{\sigma}=0$. By the parallel argument and \eqref{bb}, the assertion holds for the other elements. 
\end{proof}

Immediately, 
\begin{equation}\label{211}
P\pi^{16n+13}_{16n+29}\subset E^6\pi^{8n}_{16n+21}.
\end{equation}
Hereafter we use the following convention.

\medskip
{\bf Convention} 

In the EHP sequence arguments: 

(1)\qua Higher suspended elements in a relation are omitted. For example, in a relation 
$E^k\delta\in\{[\iota_{n-1},\beta],[\iota_{n-1},\gamma]\}$,  
if $[\iota_{n-1},\gamma]=E^l\gamma'$ for some element $\gamma'$ and $l\geq k+1$, then $[\iota_{n-1},\gamma]$ is 
omitted. 

(2)\qua Elements of order $2$ having independent Hopf invariants in a relation are omitted, if other elements are suspended. For example, in a relation $E^k\delta\ (\bmod\ \delta_1)\in\{[\iota_n,\beta]\}\ (k\geq 1)$, if $2\delta_1=0, H\delta_1\neq 0$ and $H[\iota_n,\beta]=0$, then $\delta_1$ disappears in the relation. 
\medskip

Now, we show the following:
\begin{prop}\label{nu*3}

\mbox{\em(1)}\qua
$H(E^{-3}[\iota_{8n+3},\alpha])=0$ if $\nu\alpha=0$.

\mbox{\em(2)}\qua
$H(E^{-3}[\iota_{8n+3},\beta]_{\neq 0})=\nu\beta$ for $\beta=\kappa,\nu^*, \bar{\sigma}, \bar{\kappa},\nu\bar{\kappa}$.

\mbox{\em(3)}
$H(E^{-3}[\iota_{8n+3},\nu\kappa]_{\neq 0})=4\bar{\kappa}$ if $\sharp[\iota_{8n},\bar{\kappa}]=8$.
\end{prop}
\begin{proof}
By \fullref{ex1}(2), it suffices to prove the non-triviality in (2) and (3). We show it for $\nu^*$. Let $n\equiv 3\ (8)$. By \fullref{tau12}(2), $[\iota_n,\nu^*]=E^3(\bar{\tau}_{n-3}\nu^*)$. ${\mathcal ASM}[\nu^*]$ and \eqref{4n2zeta} for $\bar{\zeta}$ induce $E^2(\bar{\tau}_{n-3}\nu^*)\in\{[\iota_{n-1},\bar{\sigma}]\}\subset E^3\pi^{n-4}_{2n+13}$ (\fullref{T11}(1)),
$E(\bar{\tau}_{n-3}\nu^*)\in P\pi^{2n-3}_{2n+17} 
=\{E(\tau_{n-3}\bar{\kappa})\}$,
$\bar{\tau}_{n-3}\nu^*\ (\bmod\ \tau_{n-3}\bar{\kappa})\in P\pi^{2n-5}_{2n+16}$ and hence, ${\mathcal CDR}[\sigma^3\ (\bmod\ \eta\bar{\kappa})=0]$. By the parallel argument, (2) for the other elements follows. We show (3). Assume that $E^3(\bar{\tau}_{8n}\nu\kappa)=[\iota_{8n+3},\nu\kappa]=0$. Then, $E^2(\bar{\tau}_{8n}\nu\kappa)\in\{[\iota_{8n+2},4\nu^*],[\iota_{8n+2},\eta\bar{\mu}]\}=0$ and\\  $E(\bar{\tau}_{8n}\nu\kappa)\in\{[\iota_{8n+1},\bar{\zeta}],[\iota_{8n+1},\bar{\sigma}]\}=\{E(\tau_{8n}\bar{\zeta}),E(\tau_{8n}\bar{\sigma})\}$. This and the assumption $\sharp[\iota_{8n},\bar{\kappa}]=8$ imply $\bar{\tau}_{8n}\nu\kappa+a\tau_{8n}\bar{\zeta}+b\tau_{8n}\bar{\sigma}\in\{4[\iota_{8n},\bar{\kappa}]\}$. Since 
$\eta\bar{\zeta}=\eta\bar{\sigma}=0$, we get 
${\mathcal CDR}[\nu^2\kappa=0]$.
\end{proof} 

Immediately, 
$$
[\iota_{8n+7},\sigma^3]=0.
$$
By \fullref{nT}(3), we have:
\begin{lem}\label{nT11}
$H(E^{-3}[\iota_{4n+2},\alpha])=\frac{1+(-1)^n}{2}\nu\alpha$ 
for $\alpha=\bar{\nu},\varepsilon,\kappa,\nu\kappa, 
\nu^2\kappa,\bar{\sigma},\nu\bar{\sigma}$. In particular, $H(E^{-3}[\iota_{8n+2},\alpha])=\nu\alpha$ 
for $\alpha=\kappa,\nu\kappa,\nu^2\kappa,\bar{\sigma}$.
\end{lem}

Immediately,
\begin{gather}\label{nk6}
[\iota_{8n+6},\nu\kappa]=0\\
[\iota_{8n+6},\nu\bar{\sigma}]=0.\tag*{\hbox{and}}
\end{gather}
We need the following:
\begin{lem}\label{nT22}\

\mbox{\em (1)}\qua
$H(E^{-3}[\iota_{4n+1},\alpha])=\frac{1+(-1)^n}{2}\nu\alpha$ for  $\alpha=\nu,\nu^2,\nu\kappa,\nu^*,\bar{\sigma},\nu^2\kappa$,\\ 
$\nu\bar{\sigma},\nu\bar{\kappa}$ and  
$H(E^{-3}[\iota_{4n+1},\beta_1])=
H(E^{-3}[\iota_{4n+1},\eta\beta_2])=0$ for $\beta_1=\zeta,\bar{\zeta},\zeta_{3,\ast}$;\\  $\beta_2=\mu,\bar{\mu},\mu_{3,\ast}$. In particular, $H(E^{-3}[\iota_{8n+1},\alpha])=\nu\alpha$\ for $\alpha=\nu,\nu^2,\nu\kappa,\nu^*,\bar{\sigma}$,\\ 
$\nu^2\kappa,\nu\bar{\kappa}$.

\mbox{\em (2)}\qua
$H(E^{-4}[\iota_{8n+5},\delta_1])=0$ for $\delta_1=\eta^2,\nu,\eta^2\sigma,\eta\varepsilon,\eta^2\rho,\nu\kappa,\eta\mu,\eta\eta^*,\eta\bar{\mu}$. 

\mbox{\em (3)}\qua
$H(E^{-4}[\iota_{8n+1},\delta_2])=0$ for $\delta_2=\nu^3,\eta^2\sigma,\sigma^2,\eta^2\rho,\nu\bar{\sigma},\eta\eta^*\sigma,\eta^2\bar{\rho},\nu^2\bar{\kappa}$. 

\mbox{\em (4)}\qua
$H(E^{-6}[\iota_{8n+5},\eta^2\delta_3])=0$ for  $\delta_3=\rho,\bar{\rho}$. 
\end{lem}
\begin{proof}
(1) is a direct consequence of  \fullref{nT}(2). Let $n\equiv 5\ (8)$. Then, $\P^{n-1}_{n-5}=E^{n-5}\P^4_0$ and we can take $\lambda_{n-5,5}=\tilde{\imath}\eta$. 
By the relations $\eta\delta_1=0,4\tilde{\imath}=i\eta^2$ \eqref{2tilio} and 
$\eta\delta=0\ (\delta=\nu,\zeta,\nu^*,\bar{\zeta})$, we have $\lambda_{n-5,5}\circ\delta_1=0$. Hence, \fullref{nTkl}[n-8;5,0] leads to (2). 

In Proposition [n-5;5,4], $\P^{n-1}_{n-5}=E^{n-9}\P^8_4$ for $n\equiv 1\ (8)$. By $(\mathcal{P}^8_4)$, we have 
$$
\pi^s_8(\P^8_4)=i''_*\pi^s_8(\P^7_4)=\{i''s_2\eta,i''ti\nu\}\cong(\Z_2)^2\ (i=i^{4,7}_4,i''=i^8_4).\leqno{(\ast)}
$$ 
So, we take 
\begin{equation}\label{gam84}
\gamma_{8,4}=i''(s_2\eta+ti\nu)
\end{equation}  
and 
$\lambda_{n-5,5}\circ\delta_2=0$. 

In \fullref{nTkl}[n-7;7,6], $\P^{n-1}_{n-7}=E^{n-13}\P^{12}_6$ for $n\equiv 5\ (8)$. Since $\P^{12}_6/\P^7_6=\P^{12}_8=E^8\P^4_0$, we have ${p^{12}_{8,6}}_*(\lambda_{n-7,7}\circ\eta^2)\in\pi^s_4(\P^4)\circ\eta^2=0$ and
$\lambda_{n-7,7}\circ\eta^2\in{i^{7,12}_6}_*\pi^s_{14}(\P^7_6)$. Hence, by the fact that $\pi^s_{14}(\P^7_6)\cong
\pi^s_8\oplus\pi^s_7$ and $\pi^s_8\circ\delta_3=\pi^s_7\circ\delta_3=0$, we obtain $\lambda_{n-7,7}\circ\eta^2\delta_3=0$. 
\end{proof}

By \fullref{nT22}(1), we obtain 
\begin{equation}\label{sgm35}
[\iota_{8n+5},\sigma^3]=0,
\end{equation}
$$
[\iota_{8n+5},\nu\bar{\sigma}]
=0
$$
and
\begin{equation}\label{star2}
P\pi^{8n+3}_{8n+21+k}\subset E^3\pi^{4n-2}_{8n+16+k}\ (k=0,1).  
\end{equation}
We also note the following. 

\begin{rem}
$H(E^{-3}[\iota_{8n+1},\nu\kappa])=\nu^2\kappa$,\ while $[\iota_{8n+1},\nu\kappa]=0$\ \eqref{nukap5}. 
\end{rem}

\section[Concerning Nomura's results]{Concerning Nomura's results \cite{N}}

In this section, we recollect Nomura's results \cite{N}, prove a part of them by using \fullref{nTkl} and add results needed in the next section. By use of the cell structures of $\P^{n-1}_{n-k}$, we determine some group structures of $\pi^s_{n-1}(\P^{n-1}_{n-k})$ for $4\leq k\leq 8$, which overlap with \cite[Section 3]{Saito}. First we show the result including the known one \cite[4.10;18]{N}. 

\begin{lem} \label{nT1M}
$H(E^{-7}[\iota_{16n+3},\sigma])=\sigma^2$ and 
$H(E^{-7}[\iota_{16n+k},\sigma^2])=\sigma^3$ for  $k=0,1,3,7$.
\end{lem}
\begin{proof}
Let $n\equiv 0\ (16)$. By \eqref{des7}, $[\iota_n,\sigma^2]
=\sigma_n\circ[\iota_{n+7},\iota]
=E^7(\sigma_{n-7}\delta_n)$ and 
$H(\sigma_{n-7}\delta_n)$ $=\sigma^3_{2n-15}$. Let $n\equiv 7\ (16)$. By \eqref{des7}, $[\iota_n,\sigma^2]
=E^7(\delta_{n-7}\sigma)$ and 
$H(\delta_{n-7}\sigma)=\sigma^3_{2n-15}$. Let $n\equiv 1\ (16)$. We have $\P^{n-1}_{n-8}=E^{n-17}\P^{16}_9$ and $\P^{n-2}_{n-8}=E^{n-17}\P^{15}_9=E^{n-9}\P^7=E^{n-9}\P^6\vee S^{n-2}$. By inspecting \cite[Proposition 4.3]{Mu2}, $$\pi^s_8(\P^6)=\{\tilde{\eta}'\eta,\widetilde{i\nu}, i^{2,6}\widetilde{\nu^2}, i\sigma\}\cong(\Z_2)^4,
$$ 
where $\widetilde{i\nu}\in\langle i^{4,6}\tilde{\imath},\bar{\eta},i\nu
\rangle$ and   
$\widetilde{i\nu}\circ\sigma\in\langle i^{4,6}\tilde{\imath},\bar{\eta},i\nu\rangle\circ\sigma
=i^{4,6}\tilde{\imath}\circ\langle\bar{\eta},i\nu,\sigma\rangle=0.$
So, by \fullref{coexnu2}(8), 
$\pi^s_8(\P^6)\circ\sigma^2=\{i\sigma^3\}$. 
Since $p_*\co  \pi^s_{16}(\P^{16}_9)\to\pi^s_{16}(S^{16})$ is trivial, $\pi^s_{16}(\P^{16}_9)=i_*\pi^s_{16}(\P^{15}_9)$.
This implies  $(\pi^s_{16}(\P^{16}_9)-\{i\sigma\})\circ\sigma^2=0$ and 
hence, by \fullref{nT3}[[n-8;8,9]], the assertion follows. 

Next, let $n\equiv 3\ (16)$. In \fullref{nT3}[[n-8;8,11]],
 $\P^{n-1}_{n-8}=E^{n-11}\P^{10}_3$. Since
 $\{E\P^4,\P^2\}=\{i\bar{\eta}',\tilde{\eta}\bar{\eta}p^4_3,i\nu
 p\}\cong(\Z_2)^3$, $\P^7_3=\P^6_3\vee S^7$ and $Sq^4$ is trivial
 on\break
 $\tilde{H}^3(\P^8_3;\Z_2)$, we can take  
$\P^8_3=M^4\cup_{i\bar{\eta}'}C(E^3\P^4).$
From the relations $\bar{\eta}'\tilde{\imath}\eta=0$ and $\bar{\eta}'\tilde{\imath}\nu=\nu^2$ \eqref{etatild}, we obtain 
$\pi^s_8(\P^8_3)=\{\widetilde{\tilde{\imath}\eta},\widetilde{i\nu}'\}$ and $\pi^s_9(\P^8_3)=\{\widetilde{\tilde{\imath}\eta}\eta\}\cong\Z_2$, where 
$$
\widetilde{\tilde{\imath}\eta}\in\langle i',\bar{\eta}',\tilde{\imath}\eta\rangle \hspace{5mm} \mbox{and} \hspace{5mm}  
\widetilde{i\nu}'\in\langle i',\bar{\eta}',i\nu\rangle\ (i'=i^{3,8}_3).
$$
By \eqref{gam85}, we obtain 
$$
\gamma_{8,3}=\widetilde{\tilde{\imath}\eta}+\widetilde{i\nu}'.
$$
By the fact that $\pi^s_6(\P^4)=\{\tilde{\imath}\nu\}\cong\Z_2$ and \eqref{etatild}, we obtain $\langle\bar{\eta}',i\nu,\eta\rangle=\pi^s_6=\{\nu^2\}$ and
$\widetilde{i\nu}'\eta\in i'\circ\langle \bar{\eta}',i\nu,\eta\rangle=\{i'\nu^2\}=0$. Hence, 
$$
\gamma_{8,3}\eta=\widetilde{\tilde{\imath}\eta}\eta
$$
and $\pi^s_{10}(\P^{10}_3)=i'''_*\pi^s_{10}(\P^8_3)=i'''_*i''_*\pi^s_{10}(M^4)=\{i'''i''\widetilde{\nu^2},i\sigma\}\ (i'''=i^{8,10}_3,i''=i^{4,8}_3)$. Therefore, by \fullref{coexnu2}(8), $(\pi^s_{10}(\P^{10}_3)-\{i\sigma\})\circ\sigma=0$. This implies 
$H(E^{-7}[\iota_n,\sigma])=\sigma^2$ and  
$H(E^{-7}[\iota_n,\sigma^2])=\sigma^3$ \eqref{bb}. 
\end{proof}
 
\begin{equation} 
[\iota_{16n+11},\sigma^2]=0,\ 
[\iota_{16n+k},\sigma^3]=0\ (k=8,11,15)\tag*{\hbox{Immediately,}} 
\end{equation}and 
\begin{equation}\label{sgm39}
[\iota_{16n+9},\sigma^3]=0.
\end{equation}
Next, we show the following \cite[Table 2, 4.15;16]{N}. 
\begin{lem}\label{nTa11}\

\mbox{\em(1)}\qua 
$H(E^{-4}[\iota_{8n+4},16\rho])=\bar{\zeta}$.

\mbox{\em(2)}\qua
$H(E^{-5}[\iota_{8n+3},\nu\kappa])=\eta^2\bar{\kappa}$. 

\mbox{\em(3)}\qua
$H(E^{-6}[\iota_{8n},\nu^3])=\eta\kappa$.
\end{lem}
\begin{proof}
Let $n\equiv 4\ (8)$. 
In \fullref{nTkl}[n-5;5,7], $\P^{n-1}_{n-5}=E^{n-12}\P^{11}_7$. Let $\widetilde{2\iota}\in\langle i',i\nu,2\iota\rangle\ ( i=i^{7,10}_7,i'=i^{11}_7)$ be a coextension of $2\iota$ in 
$(\mathcal{P}^{11}_7)$. By \fullref{P76}, we can take   
\begin{equation}\label{gam117}
\gamma_{11,7}=\widetilde{2\iota}
+i^{11}_7\tilde{\eta}''. 
\end{equation}
By \fullref{coexnu2}(3), $\lambda_{n-5,5}\circ 16\rho\in i'\circ\langle i\nu,2\iota,16\rho\rangle=i\bar{\zeta}$. 

In \fullref{nTkl}[n-6;6,1], $\P^{n-1}_{n-6}=E^{n-11}\P^{10}_5$ for $n\equiv 3\ (8)$. By the cell structure  $(\mathcal{P}^4)$, we obtain $\{M^5,\P^4\}=\{\tilde{\imath}\bar{\eta},
i^{2,4}\overline{i\nu}\}\cong\Z_4\oplus\Z_2$, where $2\tilde{\imath}\bar{\eta}=\tilde{\imath}\eta^2p$. Since $Sq^k$ on $\tilde{H}^{9-k}(\P^{10}_5;\Z_2)$ is non-trivial for $k=2,4$, 
$$
(\mathcal{P}^{10}_5) \hspace{1cm} \P^{10}_5=\P^8_5\cup_{\tilde{\imath}\bar{\eta}
+i^{2,4}\overline{i\nu}}CM^9\ (\P^8_5=E^4\P^4).
$$
From the natural isomorphisms   $\pi^s_{10}(\P^{10}_5)\cong\pi^s_{10}(\P^8_5)\cong\pi^s_6(\P^4)=\{\tilde{\imath}\nu\}\cong\Z_2$, we obtain 
$$
\pi^s_{10}(\P^{10}_5)=\{i'\tilde{\imath}\nu\}\cong\Z_2\  (i'=i^{8,10}_5),
$$
\begin{equation}\label{ast2}
\gamma_{10,5}=i'\tilde{\imath}\nu 
\end{equation}  
and
$$
(\mathcal{P}^{11}_5) \hspace{1.5cm} 
\P^{11}_5=\P^{10}_5\cup_{i'\tilde{\imath}\nu}e^{11}.
$$
Hence, by the relation $4\bar{\kappa}=\nu^2\kappa$ and \eqref{2tilio}, $\lambda_{n-6,6}\circ\nu\kappa=
4i'\tilde{\imath}\bar{\kappa}=i\eta^2\bar{\kappa}$. 

 In \fullref{nTkl}[n-7;7,1], $\P^{n-1}_{n-7}=E^{n-8}\P^7$ for $n\equiv 0\ (8)$. Let $s_3\co  S^7\hookrightarrow\P^7=\P^6\vee S^7$ be the canonical inclusion. Then, we take 
\begin{equation}\label{gam7} 
\gamma_7=2s_3+i^{6,7}\tilde{\eta}'. 
\end{equation}
By \fullref{P76}(1), $\tilde{\eta}'\circ\nu^3\in i_{4,6}\circ\langle\tilde{\imath}\bar{\eta},\tilde{\eta},\nu^3\rangle.$ By \eqref{2tilio} and Lemmas \ref{coexnu2}(6),(9), \ref{P76}(3), 
$$
\langle\tilde{\imath}\bar{\eta},\tilde{\eta},\nu^3\rangle
\subset\langle\tilde{\imath},2\nu,\nu^3\rangle\supset i^{2,4}\circ\langle\tilde{\eta},\nu,\nu^3\rangle
=i\eta\kappa
$$
$$
(\bmod\ \tilde{\imath}\circ\pi^s_{13}
+\pi^s_7(\P^4)\circ\eta\bar{\nu}=0).
$$
Hence,
$\lambda_{n-7,7}\circ\nu^3=i\eta\kappa$. 
\end{proof}
 
Immediately,
$$
[\iota_{8n+7},\bar{\zeta}]
=[\iota_{8n+5},\eta^2\bar{\kappa}]
= [\iota_{8n+1},\eta\kappa]=0.
$$
By the way, the argument in \cite[Section 4]{GM} implies that 
$\Delta_\H\co  \pi_{8n+10}(S^{8n+7})\to\pi_{8n+9}(Sp(2n+1))$ is trivial on the the $2$ primary component and 
$$
\Delta(\eta^2_{8n+5}\bar{\kappa})=4i_*\Delta_\H(\bar{\kappa}_{8n+7})=i_*\Delta_\H(\nu_{8n+7})\nu\kappa=0,
$$
where $\Delta_\H$ is the the symplectic connecting map and $i\co  Sp(2n+1)\hookrightarrow SO(8n+7)$ the canonical inclusion. 

The non-triviality of $[\iota_{8n},\nu^3]$ is proved in \cite{GM}.

Now we show the following result overlapping with  \cite[4.12]{N}.  

\begin{lem} \label{Ws25}
$H(E^{-4}[\iota_{8n+4},\sigma^2])=\nu^*$ and 
$H(E^{-5}[\iota_{8n+5},\sigma^2])=\bar{\sigma}$.
\end{lem}
\begin{proof}
In \fullref{nTkl}[n-5;5,7], $\P^{n-1}_{n-5}=E^{n-12}\P^{11}_7$ for $n\equiv 4\ (8)$. By Lemmas \ref{coexnu2}(2), \ref{P76}(3) and \eqref{gam117},  
$\lambda_{n-5,5}\circ\sigma^2=i'\langle i\nu,2\iota,\sigma^2\rangle=i\nu^*$. 
 
In \fullref{nTkl}[n-6;6,7], $\P^{n-1}_{n-6}=E^{n-13}\P^{12}_7$ for $n\equiv 5\ (8)$. We see that  $\{M^7,\P^6_3\}=\{i'\overline{i\nu},i'\tilde{\eta}\bar{\eta},\tilde{\eta}''p\}\cong(\Z_2)^3\ (i'=i^{4,6}_3)$. By $(\mathcal{P}^{11}_7)$, we have
$$
(\mathcal{P}^{12}_7) \hspace{1.5cm} 
\P^{12}_7=\P^{10}_7\cup_{i'\overline{i\nu}+\tilde{\eta}''p}CM^{11}
$$
and
$\pi^s_{12}(\P^{12}_7)=\{\widetilde{i\eta},i''\widetilde{i\nu}\}\cong(\Z_2)^2$, where $\widetilde{i\eta}\in\langle i'',i'\overline{i\nu}+\tilde{\eta}''p,i\eta\rangle\ (i''=i^{10,12}_7)$ and $\widetilde{i\nu}\in\langle i',\bar{\eta},i\nu\rangle\in\pi^s_{12}(\P^{10}_7)$. 
Since $\langle i'',i'\overline{i\nu}+\tilde{\eta}''p,i\eta\rangle\supset\langle i'',(i'\overline{i\nu}+\tilde{\eta}''p)\circ i,\eta\rangle
=\langle i'',i'i\nu,\eta\rangle\supset\langle i''i'i,\nu,\eta\rangle$, we can choose $\widetilde{i\eta}$ such that
\begin{equation}\label{star}  
\widetilde{i\eta}\in\langle i''',\nu,\eta\rangle\ (i'''=i''i'i=i^{7,12}_7).
\end{equation}
From the fact that $Sq^4$ is trivial on $\tilde{H}^9(\P^{13}_7;\Z_2)$, we take 
$\gamma_{12,7}=\widetilde{i\eta}$ and 
$$
\widetilde{i\eta}\circ\sigma^2\in i'''\circ\langle\nu,\eta,\sigma^2\rangle=i'''\bar{\sigma}\ 
(\bmod\ 0).
$$
This implies $\lambda_{n-6,6}\circ\sigma^2
=i'''\bar{\sigma}$. 
\end{proof}

Immediately,
$$
[\iota_{8n+7},\nu^*]=[\iota_{8n+7},\bar{\sigma}]=0.
$$
Next, we prove the following \cite[4.13;14;16,Table 2]{N}. 
\begin{lem}\label{No} \

\mbox{\em(1)}\qua
$H(E^{-5}[\iota_{8n+2},\eta])=\nu^2$.

\mbox{\em(2)}\qua 
$H(E^{-6}[\iota_{8n+1},\eta^2])=\varepsilon$.

\mbox{\em(3)}\qua
$H(E^{-5}[\iota_{8n+2},\eta^*])=\sigma^3$. 

\mbox{\em(4)}\qua 
$H(E^{-6}[\iota_{8n+1},\eta\eta^*])=\eta^*\sigma$.

\mbox{\em(5)}\qua
$H(E^{-6}[\iota_{8n+6},\kappa])=\bar{\kappa}$. 

\end{lem}
\begin{proof} 
In \fullref{nTkl}[n-6;6,6], $\P^{n-1}_{n-6}=E^{n-10}\P^9_4$ for $n\equiv 2\ (8)$, $\pi^s_9(\P^9_4)\cong\pi^s_9(\P^9_5)\cong\Z\ (\mathcal{P}^9_5)$ and $\gamma_{9,4}\circ\eta^*\in i'''\circ\langle\gamma_{8,4},2\iota,\eta^*\rangle\ (i'''=i^9_4)$ \eqref{gam2n}. 
By the relations $\langle p,i,2\iota\rangle=\pm\iota$,  $\langle i\nu,2\iota,\eta\rangle=0$, \eqref{gam74} and \eqref{gam84}, 
we have $2i''s_2=i''(i\nu\pm t\tilde{\eta})\ (i''=i^8_4)$ 
and 
$$
\langle\gamma_{8,4},2\iota,\eta\rangle\subset \langle i''s_2\eta,2\iota,\eta\rangle+\langle i''ti\nu,2\iota,\eta\rangle
\ni\pm 2i''s_2\nu=i\nu^2
$$
$$ 
(\bmod\ i''(s_2\eta+ti\nu)
\circ\pi^s_2+\pi^s_9(\P^8_4)\circ\eta).
$$
The indeterminacy is trivial, because $\pi^s_9(\P^8_4)=\{i''s_2\eta^2\}\cong\Z_2$ and  $i''s_2\eta^3=4i''s_2\nu=0$. This implies $\lambda_{n-6,6}\circ\eta=i\nu^2$.  

In \fullref{nTkl}[n-7;7,2], $\P^{n-1}_{n-7}=E^{n-9}\P^8_2$ for $n\equiv 1\ (8)$. We obtain  $\{M^5,\P^4_2\}$\break
$=\{\tilde{\imath}'\bar{\eta},i\nu\}\cong\Z_4\oplus\Z_2$, 
$\pi^s_6(\P^6_2)=\{i'\tilde{\imath}'\nu\}\cong\Z_2\ (i'=i^{4,6}_2)$ and 
$\pi^s_7(\P^6_2)=\{\tilde{\eta}'''\}\cong\Z_8$, where 
$\tilde{\eta}'''\in\langle i',\tilde{\imath}'\bar{\eta},\tilde{\eta}\rangle$ and $2\tilde{\eta}'''\in\langle i',i\eta\bar{\eta},\tilde{\eta}\rangle\ (i=i^{2,4}_2)$. 
We also obtain $\{M^7,\P^6_2\}=\{i'\overline{\tilde{\imath}'\nu},\tilde{\eta}'''p\}\cong(\Z_2)^2$. By the cell structures
$$
\P^6_2=\P^4_2\cup_{\tilde{\imath}'\bar{\eta}}CM^5 \hspace{3mm} \mbox{and} \hspace{3mm} \P^8_2=\P^6_2\cup_{\tilde{\eta}'''p}CM^7,
$$
we have $\pi^s_8(\P^6_2)=\{\tilde{\eta}'''\eta,\widetilde{i\nu}'',i'\tilde{\imath}\nu^2\}\cong(\Z_2)^3$ and 
$\pi^s_8(\P^8_2)=\{\tilde{\imath}'''\eta\}\oplus i''_*\pi^s_8(\P^6_2)$, where $\widetilde{i\nu}''\in\langle i'\tilde{\imath}',\bar{\eta},i\nu\rangle$, 
$\tilde{\imath}'''\in\langle i'',\tilde{\eta}'''p,i\rangle\ (i''=i^{6,8}_2)$ and $2\tilde{\imath}'''=i''\tilde{\eta}'''$  \cite[Proposition 4.2]{Mu2}. We can take $\gamma_{8,2}\equiv \tilde{\imath}'''\eta\ (\bmod\ i''_*\pi^s_8(\P^6_2))$. Since $\widetilde{i\nu}''\circ\eta\in i'\tilde{\imath}'\circ\langle\bar{\eta},i\nu,\eta\rangle=i'\tilde{\imath}'\nu^2$, we obtain $\widetilde{i\nu}''\circ\eta^2=0$. By \fullref{coexnu2}(8),  $\gamma_{8,2}\circ\eta^2=2i''\tilde{\eta}'''\nu\in i''i'\circ\langle i\eta\bar{\eta},\tilde{\eta},\nu\rangle=i\varepsilon$. 

By the same argument as (1) and by \fullref{coexnu2}(4), 
$\lambda_{n-6,6}\circ\eta^*=i\sigma^3$. 
By the same argument as (2) 
and by \fullref{coexnu2}(1), 
$\gamma_{8,2}\eta\eta^*=i\eta^*\sigma$. 

Since $\langle\eta,2\iota,\kappa\rangle\ni 0$, we can choose 
a coextension $\tilde{\kappa}\in\pi^s_{16}(\P^2)$ satisfying 
$\bar{\eta}\tilde{\kappa}=0$. Notice that 
$\langle\nu,\eta,\eta\kappa\rangle
=\pm 2\bar{\kappa}$ and $\langle\nu,\bar{\eta},\tilde{\kappa}\rangle=\pm\bar{\kappa}$. In \fullref{nTkl}[n-7;7,7], $\P^{n-1}_{n-7}=E^{n-14}\P^{13}_7$ for $n\equiv 6\ (8)$. By use of $(\mathcal{P}^{12}_7)$, we get that 
$$
\pi^s_{13}(\P^{12}_7)=\{\widetilde{i\eta}\eta,i'\tilde{\eta}''\eta^2,i\nu^2\}\cong(\Z_2)^3\ (i'=i^{10,12}_7).
$$
We obtain $\widetilde{i\eta}\eta\kappa\in i'''\circ\langle\nu,\eta,\eta\kappa\rangle
=2i\bar{\kappa}=0$. By $\eqref{star}$, there exists an extension $\widetilde{i\bar{\eta}}\in\langle i''',\nu,\bar{\eta}\rangle$ of $\gamma_{12,7}=\widetilde{i\eta}$. By \eqref{gam2n}, we obtain  $\gamma_{13,7}\circ\kappa\in i^{13}_7\circ\langle\gamma_{12,7},2\iota,\kappa\rangle\ni
 i^{13}_7\widetilde{i\bar{\eta}}\tilde{\kappa}\ (\bmod\ {i^{13}_7}_*\pi^s_{13}(\P^{12}_7)\circ\kappa=0)$.
We obtain $\widetilde{i\bar{\eta}}\tilde{\kappa}\in 
i'''\circ\langle\nu,\bar{\eta},\tilde{\kappa}\rangle=i\bar{\kappa}\ (\bmod\ i''\circ\{M^6,S^0\}\circ\tilde{\kappa}=\{i''\nu^2\kappa\}=0)$ and hence, $\lambda_{n-7,7}\circ\kappa=i\bar{\kappa}$.
\end{proof}

Immediately,
$$
[\iota_{8n+4},\sigma^3]=[\iota_{8n+2},\eta^*\sigma]
=[\iota_{8n+7},\bar{\kappa}]
=0.
$$
Given an element $\alpha\in\pi_k(S^n)$, a lift $[\alpha]\in\pi_k(SO(n+1))$ of $\alpha$ is an element 
satisfying $p_{n+1}(\R)[\alpha]=\alpha$, where $p_{n+1}(\R)\co  SO(n+1)\to S^n$ is the projection.
A lift $[\alpha]$ exists if and only if $\Delta\alpha=0\in\pi_{k-1}(SO(n))$. Let $n\equiv 7\ (8)$. We know $\Delta\nu_n=0$ \cite{M1}. Note the fact that $\Delta\kappa_n=0$ \cite[Section 5]{GM} is obtained by constructing a lift of $\kappa_n$ is given by 
$$
[\kappa_n]\in\{[\nu_n],\bar{\eta},\widetilde{\bar{\nu}}\}\subset\pi_{n+14}(SO(n+1))\ (\widetilde{\bar{\nu}}:\ \mbox{a coextension of}\ \bar{\nu}).
$$
By the parallel argument, lifts of 
$\bar{\sigma}_n$ and $\bar{\kappa}_n$ are taken as follows: 
$$
[\bar{\sigma}_n]\in\{[\nu_n],\eta,\sigma^2\}\subset\pi_{n+19}(SO(n+1)); 
$$
$$
[\bar{\kappa}_n]\in\{[\nu_n],\bar{\eta},\tilde{\kappa}\}\subset\pi_{n+20}(SO(n+1)).
$$
Hence, 
$$
\Delta\bar{\sigma}_{8n+7}=\Delta\bar{\kappa}_{8n+7}=0.
$$
We need the following result overlapping with \cite[4.14]{N}.
\begin{lem}\label{nnu*}\

\mbox{\em(1)}\qua
$H(E^{-6}[\iota_{8n+3},\alpha])=0$ if $\nu\alpha=0$.  

\mbox{\em(2)}\qua
$H(E^{-6}[\iota_{8n+4k},4\nu^*])=\eta\eta^*\sigma$ or $0$ according as $k=0$ or $1$. 
\end{lem}
\begin{proof} 
In \fullref{nTkl}[n-7;7,4], $\P^{n-1}_{n-7}=E^{n-11}\P^{10}_4$ for $n\equiv 3\ (8)$. 
We have $\{\P^4,S^1\}=\{\eta\bar{\eta}p^4_3,\nu p\}\cong(\Z_2)^2\ (p=p^4_4)$, 
$\eta\bar{\eta}p^4_3\circ(\tilde{\imath}\bar{\eta}+i^{2,4}\overline{i\nu})=\eta^2\bar{\eta}$ and $p\circ(\tilde{\imath}\bar{\eta}+(i^{2,4})\overline{i\nu})=0$. So, by the fact that $\{M^5,S^0\}=0$ and $(\mathcal{P}^{10}_5)$, $p$ is extendible on $\bar{p}\in\{\P^{10}_5,S^8\}$ and $\{\P^{10}_5,S^5\}=\{\nu\bar{p}\}\cong\Z_2$. Hence,
$$
E\P^{10}_4=S^5\cup_{\nu\bar{p}}C\P^{10}_5. 
$$
Since $(\tilde{\imath}\bar{\eta}+(i^{2,4})\overline{i\nu})\circ i\nu=i\nu^2$, we have $i'i\nu^2=0$ in $\pi^s_{11}(\P^{10}_5)\ (i'=i^{8,10}_5)$. By \fullref{coexnu2}(5), $\langle i^{2,4}\overline{i\nu},\tilde{\eta},4\iota\rangle\supset
i^{2,4}\circ\langle\overline{i\nu},\tilde{\eta},4\iota\rangle=\{i\nu^2\}.$
So, by $(\mathcal{P}^{10}_5)$ and \fullref{P76}(1),  $\pi^s_{11}(\P^{10}_5)=\{\tilde{\eta}^{IV}\}\cong\Z_8$, where 
$\tilde{\eta}^{IV}\in\langle i',\tilde{\imath}\bar{\eta}
+i^{2,4}\overline{i\nu},\tilde{\eta}\rangle$ and $4\tilde{\eta}^{IV}=i'\widetilde{\tilde{\eta}\eta^2}$. 
By the fact that $\langle p',i\bar{\eta},\tilde{\eta}\rangle=\pm\nu\ (p'=p^2_2)$ and 
$$
\langle p,i^{2,4}\overline{i\nu},\tilde{\eta}\rangle
\subset\langle p',0,\tilde{\eta}\rangle\ni 0\ (\bmod\ p'\circ\pi^s_5(\P^2)+\{\P^2,S^0\}\circ\tilde{\eta}=\{2\nu\}),
$$
we obtain 
$\bar{p}\circ\tilde{\eta}^{IV}=\pm\nu$. So, by \eqref{ast2} and the relation $\bar{p}\circ i'\tilde{\imath}=0\ (i'\tilde{\imath}\in\pi^s_7(\P^{10}_5))$, 
we conclude that  $\pi^s_{10}(\P^{10}_4)
=\{\widetilde{i'\tilde{\imath}}\nu\}\cong\Z_2$ and $\gamma_{10,4}=\widetilde{i'\tilde{\imath}}\nu$, where $\widetilde{i'\tilde{\imath}}\in\pi^s_7(\P^{10}_4)$ is a coextension of $i'\tilde{\imath}$. This leads to (1).

In \fullref{nTkl}[n-7;7,1], $\P^{n-1}_{n-7}=E^{n-8}\P^7$ for $n\equiv 0\ (8)$. By \eqref{gam7} and \fullref{coexnu2}(10), $\lambda_{n-7,7}\circ 4\nu^*=i^{4,7}\widetilde{\tilde{\eta}\eta^2}\nu^*=i^{2,7}\circ\langle\tilde{\eta}p,\tilde{\eta}\eta^2,\nu^*\rangle=
i\eta\eta^*\sigma$. Hence, $\lambda_{n-7,7}\circ 4\nu^*=i\eta\eta^*\sigma$.

In \fullref{nTkl}[n-7;7,5], $\P^{n-1}_{n-7}=E^{n-12}\P^{11}_5$ for $n\equiv 4\ (\bmod\ 8)$. By use of 
$(\mathcal{P}^{11}_5)$, we can take  
\begin{equation}\label{IV}
\gamma_{11,5}=i''\tilde{\eta}^{IV}+\widetilde{2\iota},\ \mbox{where}\  \widetilde{2\iota}\in\langle i''',\tilde{\imath}\nu,2\iota\rangle\ (i''=i^{11}_5,i'''=i^{8,11}_5).
\end{equation} 
By \fullref{coexnu2}(10), $\tilde{\eta}^{IV}\circ 4\nu^*=i'\widetilde{\tilde{\eta}\eta^2}\circ\nu^*
=i\eta\eta^*\sigma$. By the relation  
$\widetilde{2\iota}\circ\eta
\in i'''\circ\langle \tilde{\imath}\nu,2\iota,\eta\rangle$ and Lemmas \ref{coexnu2}(8),\ref{P76}(3), 
$$
\langle\tilde{\imath}\nu,2\iota,\eta\rangle
\subset\langle\tilde{\imath},2\nu,\eta\rangle
\supset\langle i'\tilde{\eta},\nu,\eta\rangle\ni i^{2,4}\widetilde{\nu^2}\ 
(\bmod\ \pi^s_7(\P^4)\circ\eta=0).
$$
Hence, 
$\widetilde{2\iota}\circ\eta=i'''i^{2,4}\widetilde{\nu^2}$ and 
$\widetilde{2\iota}\circ 4\nu^*=i'''i^{2,4}\widetilde{\nu^2}\eta\eta^*=i\eta\eta^*\sigma$ by \fullref{coexnu2}(10). 
Thus, by \eqref{IV}, $\lambda_{n-7,7}\circ 4\nu^*=0$. 
This leads to (2). 
\end{proof}

Immediately,
$$
[\iota_{8n+1},\eta\eta^*\sigma]=0.
$$
Finally, we need the following \cite[4.8;9;10;11;16;17;18]{N}.
\begin{lem}\label{Nnnu*}\

\mbox{\em(1)}\qua
$H(E^{-6}[\iota_{8n+5},\eta\kappa])=\eta\bar{\kappa}$.

\mbox{\em(2)}\qua
$H(E^{-6}[\iota_{8n+4},\nu\kappa])
=\nu\bar{\kappa}$.

\mbox{\em(3)}\qua
$H(E^{-6}[\iota_{8n+2},4\bar{\kappa}])
=\nu^2\bar{\kappa}$.

\mbox{\em(4)}\qua
$H(E^{-7}[\iota_{16n+14},\eta^*])=\eta^*\sigma\ \mbox{and}\ 
H(E^{-7}[\iota_{16n+13},\eta\eta^*])=\eta\eta^*\sigma.$

\mbox{\em(5)}\qua
$H(E^{-11}([\iota_{16n+5},\nu])=\sigma^2$.

\mbox{\em(6)}\qua
$H(E^{-13}[\iota_{16n+3},\nu^2])
=\bar{\sigma}$, $H(E^{-11}[\iota_{16n+2},\eta\sigma])
=\bar{\sigma}$ and\\ 
$H(E^{-13}[\iota_{16n+1},\nu^3])=\nu\bar{\sigma}$. 
\end{lem}

Immediately,
$$
[\iota_{8n+6},\eta\bar{\kappa}]=[\iota_{8n+5},\nu\bar{\kappa}]=[\iota_{8n+3},\nu^2\bar{\kappa}]=[\iota_{16n+9},\sigma^2]=0,
$$
$$
[\iota_{16n+5},\bar{\sigma}]=
[\iota_{16n+6},\bar{\sigma}]=[\iota_{16n+3},
\nu\bar{\sigma}]=[\iota_{16n+6},
\eta^*\sigma]=[\iota_{16n+5},\eta\eta^*\sigma]=0. 
$$

\section[Completion of the proof of \ref{main}]{Completion of the proof of
\fullref{main}}

First we show: 
\begin{prop} \label{Ws2116}
$[\iota_n,\sigma^2]\neq 0$ for $n\equiv 4,5\ (8)$ or $ 
n\equiv 0,1,3\ (16)$.
\end{prop}
\begin{proof} 
By \fullref{nT1M}, we can set $[\iota_n,\sigma^2]=E^7\delta$ and $H\delta=\sigma^3$.
Let $n\equiv 0\ (16)$. By \cite{Ad}, $[\iota_{n-1},\iota]$ desuspends seven dimensions. So, ${\mathcal ASM}[\sigma^2]$ implies $E^5\delta\in P\pi^{2n-3}_{2n+13}\subset E^6\pi^{n-8}_{2n+5}$ \eqref{211} and $E^4\delta\in P\pi^{2n-5}_{2n+12}$. By \fullref{nT22}(2), $[\iota_{n-3},\alpha]$ for $\alpha=\eta\eta^*,\eta^2\rho,\nu\kappa$ desuspends five dimensions. Hence, by the relation $H(E^{-1} [\iota_{n-3},\bar{\mu}])=\eta\bar{\mu}$, we have $E^3\delta\in\{[\iota_{n-4},4\nu^*], 
[\iota_{n-4},\eta\bar{\mu}]\}.$ 
By \fullref{nTa1}(1), $[\iota_{n-4},4\nu^*]$ desuspends four dimensions. Therefore, by the relation $H(E^{-1} [\iota_{n-4},\eta\bar{\mu}])=4\bar{\zeta}$ (\fullref{T10}(2)), $E^2\delta\in\{[\iota_{n-5},\bar{\zeta}]$, 
$[\iota_{n-5},\bar{\sigma}]\}\subset E^3\pi^{n-8}_{2n+5}$   (\fullref{nu*3}(2)). 
Hence, $E\delta
\in\{[\iota_{n-6},4\bar{\kappa}]\}\subset E^2\pi^{n-8}_{2n+5}$ (\fullref{T11}(1)), $\delta\in P\pi^{2n-13}_{2n+8}$ and ${\mathcal CDR}[\sigma^3=0]$. 

Let $n\equiv 1\ (16)$. ${\mathcal ASM}[\sigma^2]$ implies  
$E^6\delta\in\{[\iota_{n-1},\eta\kappa],[\iota_{n-1},16\rho]\}$. By \fullref{nTa1}(2), $[\iota_{n-1},16\rho]$ desuspends eight dimensions and $E^5\delta\ 
(\bmod\ E\beta)\in P\pi^{2n-3}_{2n+13}$
$=0$ for $\beta=E^{-2}[\iota_{n-1},\eta\kappa]$. So, by \eqref{diamond} 
$E^4\delta\in P\pi^{2n-5}_{2n+12}\subset E^6\pi^{n-9}_{2n+4}$ (\fullref{6bamu}), 
$E^3\delta\in P\pi^{2n-7}_{2n+11}\subset E^4\pi^{n-8}_{2n+5}$ (\fullref{nT22}(1)), 
$E^2\delta\in\{[\iota_{n-5},4\bar{\zeta}],[\iota_{n-5},\bar{\sigma}]\}$
$\subset E^3\pi^{n-8}_{2n+5}$ (\fullref{T11}(3)), 
$E\delta\in\{[\iota_{n-6},\bar{\kappa}]\}\subset E^3\pi^{n-9}_{2n+4}$ (\fullref{nu*3}(2)) and hence, ${\mathcal CDR}[\delta\in P\pi^{2n-13}_{2n+8}]$.  

Let $n\equiv 3\ (16)$. ${\mathcal ASM}[\sigma^2]$ implies  $E^6\delta\in\{[\iota_{n-1},16\rho]\}$. Since 
$H(E^{-2}[\iota_{n-1},16\rho])$ $=\bar{\mu}$ by \fullref{T11}(1), $E^5\delta\ (\bmod\ E\beta_1)\in P\pi^{2n-3}_{2n+13}=$
$\{E(\tau_{n-3}\eta\rho),E(\tau_{n-3}\eta^*)\}$ for $\beta_1=E^{-2}[\iota_{n-1},16\rho]$. So, $E^4\delta\in\{[\iota_{n-3},\alpha]\}$ 
for $\alpha\in\pi^s_{17}$. 
We obtain\break $H(E^{-1}[\iota_{n-3},\bar{\mu}])=\eta\bar{\mu}$ (\fullref{T10}(1)),  
$H(E^{-1}[\iota_{n-3},\eta\eta^*])=\eta^2\eta^*$ (\fullref{T10}(2)),
$H(E^{-2}[\iota_{n-3},\eta^2\rho])=x\bar{\zeta}\ (x: \mbox{odd})$
(\fullref{T11}(2)) and $H(E^{-3}[\iota_{n-3},\nu\kappa])$ $=\nu^2\kappa$ (\fullref{nTa1}(1)). 
This induces $E^3\delta\ (\bmod\ E\beta_2, 
E^2\delta_1)\in P\pi^{2n-7}_{2n+11}=0$ for $\beta_2
=E^{-2}[\iota_{n-3},\eta^2\rho]$ and
$\delta_1=E^{-3}[\iota_{n-3},\nu\kappa]$. 
Hence, by \eqref{4n2zeta} for $\bar{\zeta}$, $E^2\delta\ (\bmod\ E\delta_1)\in\{[\iota_{n-5},\bar{\sigma}]\}\subset E^6\pi^{n-11}_{2n+2}$ (\fullref{6bamu}),  $E\delta\in\{E(\tau_{n-7}\bar{\kappa})\}$ and   
${\mathcal CDR}[\delta\ (\bmod\ \tau_{n-7}\bar{\kappa})\in P\pi^{2n-13}_{2n+8}]$. 

Let $n\equiv 4\ (8)$. \fullref{Ws25} and ${\mathcal ASM}[\sigma^2]$ imply $E^3\delta\in P\pi^{2n-1}_{2n+14}=\{[\iota_{n-1},\rho]\}\subset E^4\pi^{n-5}_{2n+8}$ (\fullref{nu*3}(1)), for $\delta=\delta(\nu^*)=E^{-4}[\iota_n,\sigma^2]$. By \fullref{T11}(1), $H(E^{-2}[\iota_{n-2},\eta^*])=2\nu^*$ and 
$[\iota_{n-2},\eta\rho]$ desuspends three dimensions. This induces  
$E\delta\ (\bmod\ E\delta_1)\in\{[\iota_{n-3},\alpha]\}$, where $\delta_1=\delta(2\nu^*)=E^{-2}[\iota_{n-2},\eta^*]$ and\break
$\alpha=\eta^2\rho,
\eta\eta^*,\nu\kappa,\bar{\mu}$. Hence, $\delta\ (\bmod\ \delta_1,\tau_{n-4}\alpha)\in\{[\iota_{n-4},\nu^*], [\iota_{n-4},\eta\bar{\mu}]\}$ and\\ ${\mathcal CDR}[\nu^*\ (\bmod\ \eta\bar{\mu})\in \{2\nu^*\}]$. 

Let $n\equiv 5\ (8)$. \fullref{Ws25} and ${\mathcal ASM}[\sigma^2]$ induce
$E^4\delta\in\{[\iota_{n-1},\eta\kappa],[\iota_{n-1},16\rho]\}$, where $\delta=\delta(\bar{\sigma})=E^{-5}[\iota_n,\sigma^2]$. By \eqref{diamond} and  
\fullref{nTa11}(1),    
$E^3\delta\ (\bmod\ E\delta_1,E^3\delta_2)\in 
P\pi^{2n-3}_{2n+13}=0$ and    
$E^2\delta\ (\bmod\ E^2\delta_2)\in\{[\iota_{n-3},\nu\kappa],[\iota_{n-3},\bar{\mu}]\},$ where $\delta_1=\delta(\nu\kappa)=E^{-2}[\iota_{n-1},\eta\kappa]$ and $\delta_2=\delta(\bar{\zeta})=E^{-4}[\iota_{n-1},16\rho]$. By \fullref{T11}(1), $[\iota_{n-3},\nu\kappa]$ desuspends three dimensions and $H(E^{-2}[\iota_{n-3},\bar{\mu}])=2\bar{\zeta}$. Hence, 
for $\delta_3=\delta(2\bar{\zeta})=E^{-2}[\iota_{n-3},\bar{\mu}]$, we have $E\delta\ (\bmod\ E\delta_2, 
E\delta_3)\in P\pi^{2n-7}_{2n+11}=\{E(\tau_{n-5}\nu^*), E(\tau_{n-5}\eta\bar{\mu})\}$,   
$\delta\ (\bmod\ \delta_2, \delta_3, \tau_{n-5}\nu^*,\tau_{n-5}\eta\bar{\mu})\in P\pi^{2n-9}_{2n+10}$ and  ${\mathcal CDR}[\bar{\sigma}\ (\bmod\ \bar{\zeta})\in\{2\bar{\zeta}\}]$. 
\end{proof}

Next we show the following:
\begin{prop}\label{bakap8}
$H(E^{-3}[\iota_{8n},\nu^2\kappa]_{\neq 0})=
4\nu\bar{\kappa}$ and 
$H(E^{-3}[\iota_{8n+2},\nu\kappa]_{\neq 0})=4\bar{\kappa}$.
\end{prop} 
\begin{proof}
Let $n\equiv 0\ (8).$ 
By \fullref{nTa1}(1) and \eqref{bb}, $H(E^{-3}[\iota_n,\nu^2\kappa])=\nu^3\kappa=
4\nu\bar{\kappa}\ (\delta\kappa=E^{-3}[\iota_n,\nu^2\kappa])$ for $\delta=E^{-3}[\iota_n,\nu^2]$. Then, 
 ${\mathcal ASM}[\nu^2\kappa]$ induces
 $E^2(\delta\kappa)\in P\pi^{2n-1}_{2n+20}=0$, $E(\delta\kappa)\in P\pi^{2n-3}_{2n+19}=\{[\iota_{n-2},\nu\bar{\sigma}]\}\subset E^6\pi^{n-8}_{2n+11}$ (\fullref{6bamu}), and hence $\delta\kappa\in\ P\pi^{2n-5}_{2n+18}$ and ${\mathcal CDR}[4\nu\bar{\kappa}=0]$.

Next, let $n\equiv 2\ (8)$. 
By \fullref{nT11}, there exists an element 
$\delta\in\pi^{n-3}_{2n+13}$ such that $[\iota_n,\nu\kappa]=E^3\delta$ and 
$H\delta=\nu^2\kappa$. 
Hence, 
${\mathcal ASM}[\nu\kappa]$ and \eqref{star2} induce 
$E\delta\in\{[\iota_{n-2},4\bar{\zeta}], [\iota_{n-2},\bar{\sigma}]\}\subset E^2\pi^{n-4}_{2n+12}
$ (\fullref{T11}(3)) and ${\mathcal CDR}[\delta\in P\pi^{2n-5}_{2n+15}=0]$. 
\end{proof}

By Propositions \ref{nu*3}(3), \ref{bakap8} and 
the properties of Whitehead products,  
$$
\sharp[\iota_{8n},\bar{\kappa}]=8 \hspace{5mm} \mbox{and} \hspace{5mm}  
\sharp[\iota_{8n},\nu\kappa]=
\sharp[\iota_{8n+3},\nu\kappa]=\sharp[\iota_{8n+2},\kappa]=2.  
$$
We show: 
\begin{prop}\label{kap}
$\sharp[\iota_{8n+6},\kappa]=
\sharp[\iota_{8n+5},\eta\kappa]=\sharp[\iota_{8n+4},\nu\kappa]=2$. 
\end{prop}
\begin{proof}
Let $n\equiv 6\ (8)$. \fullref{No}(5) and 
${\mathcal ASM}[\kappa]$ imply $E^5\delta\in\{[\iota_{n-1},\rho],[\iota_{n-1},\eta\kappa]\}$ for $\delta=\delta(\bar{\kappa})=E^{-6}[\iota_n,\kappa]$. 
By the relation $H(\tau_{n-2}\rho)$
$=\eta\rho$ and \fullref{Nnnu*}(1), $E^4\delta\in\{[\iota_{n-2},\eta\rho],[\iota_{n-2},\eta^*]\}$. By \fullref{T10}(1), 
\begin{equation}\tag{$\star$}
H(E^{-1}[\iota_{n-2},\eta\rho])=\eta^2\rho;\ 
H(E^{-1}[\iota_{n-2},\eta^*])=\eta\eta^*.\label{starX}
\end{equation}
Therefore, $E^3\delta\in
P\pi^{2n-5}_{2n+12}=\{E^3(\bar{\tau}_{n-6}\nu\kappa)\}$. By the fact
that $[\iota_{n-4},\eta\bar{\mu}]=[\iota_{n-4},\eta^2\eta^*]$ $=0$ and \eqref{star2}, $E^2\delta\  (\bmod\ E^2(\bar{\tau}_{n-6}\nu\kappa))=0$, $E\delta\ (\bmod\ E(\bar{\tau}_{n-6}\nu\kappa))\in P\pi^{2n-9}_{2n+10}\subset E^3\pi^{n-8}_{2n+5}$, $\delta\ (\bmod\ \bar{\tau}_{n-6}\nu\kappa)
\in\{[\iota_{n-6},\bar{\kappa}]\}$ and 
hence, ${\mathcal CDR}[\bar{\kappa}\in \{2\bar{\kappa}\}]$. 

Let $n\equiv 5\ (8)$. 
\fullref{Nnnu*}(1) and ${\mathcal ASM}[\eta\kappa]$ imply  
$E^5\delta\in\{[\iota_{n-1},\eta\rho], [\iota_{n-1},\eta^*]\}$ for $\delta=\delta(\eta\bar{\kappa})=E^{-6}[\iota_n,\eta\kappa]$. 
By \eqref{starX}, $E^4\delta\in P\pi^{2n-3}_{2n+14}=\{[\iota_{n-2},\nu\kappa]\}\subset E^5\pi^{n-7}_{2n+7}$ (\fullref{nTa11}(2)), 
$E^3\delta\in\{[\iota_{n-3}, 4\nu^*],$ 
$[\iota_{n-3}, \eta\bar{\mu}]\}=0$, 
$E^2\delta\in P\pi^{2n-7}_{2n+12}\subset E^3\pi^{n-7}_{2n+7}$ \eqref{star2} and  
$E\delta\in\{[\iota_{n-5},4\bar{\kappa}]\}
\subset E^3\pi^{n-10}_{2n+9}$ (\fullref{T11}(3)). 
Hence, ${\mathcal CDR}[\delta\in P\pi^{2n-11}_{2n+10}=0]$. 

Let $n\equiv 4\ (8)$.  $E^5\delta\in\{E^3(\bar{\tau}_{n-4}\nu^*)\}$
for $\delta=\delta(\nu\bar{\kappa})=E^{-6}[\iota_n,\nu\kappa]$. By the
relation $H(E^{-3}[\iota_{n-2},\bar{\sigma}])=\nu\bar{\sigma}$
(\fullref{nT11}) and \eqref{4n2zeta} for $\bar{\zeta}$, $E^4\delta\
(\bmod\ E^2(\bar{\tau}_{n-4}\nu^*))\in\{E^3\delta_1\}$ and $E^3\delta\
(\bmod\
E(\bar{\tau}_{n-4}\nu^*),E^2\delta_1)\in\{E(\tau_{n-4}\bar{\kappa})\}$,
where $\delta_1=\delta(\nu\bar{\sigma})=$\break
$E^{-3}[\iota_{n-2},\bar{\sigma}]$. From the relations
$H(\bar{\tau}_{n-4}\nu^*)=\sigma^3$,
$H(\tau_{n-4}\bar{\kappa})=\eta\bar{\kappa}$ and\break
$H(E^{-1}[\iota_{n-4},\eta\bar{\kappa}])=\eta^2\bar{\kappa}$
(\fullref{T10}(1)), we obtain $E^2\delta\ (\bmod\ E\delta_1)\in$\break
$\{[\iota_{n-4},\sigma^3]\}\subset E^7\pi^{n-11}_{2n+5}$
(\fullref{nTa1}(3)), $E\delta\in P\pi^{2n-9}_{2n+13}=0$, $\delta\in
P\pi^{2n-11}_{2n+12}$ and hence, ${\mathcal CDR}[\nu\bar{\kappa}\in
2\pi^s_{23}]$.
\end{proof}

Since $[\iota_{8n+4},\nu^2]=0$,  $[\iota_{8n+6},\nu\kappa]=0$ \eqref{nk6} and 
$H[\iota_{2n},\bar{\kappa}]=\pm 2\bar{\kappa}$, we have
\begin{gather*}
\sharp[\iota_{8n+k},\bar{\kappa}]=4\ \mbox{for}\ k=4,6.\\
\sharp[\iota_{8n+4},\nu\bar{\kappa}]=4.\tag*{\hbox{Similarly,}}
\end{gather*}
Now, we show:
\begin{prop}\label{nu**}
$\sharp[\iota_{8n+2},\eta^*]=\sharp[\iota_{8n+1},\nu^*]=\sharp[\iota_{8n},4\nu^*]=2$. 
\end{prop}
\begin{proof}
Let  $n\equiv 2\ (8)$. By \eqref{nukap5} and 
\fullref{No}(2);(3), $[\iota_{n-1},\alpha]\in E^6\pi^{n-7}_{2n+8}$ for $\alpha=\nu\kappa, \eta^2\rho$ and $\eta\eta^*$.  So, ${\mathcal ASM}[\eta^*]$ induces $E^4\delta\in\{E(\tau_{n-2}\bar{\mu})\}$ and   
$E^3\delta\in\{[\iota_{n-2},4\nu^*], [\iota_{n-2},\eta\bar{\mu}]\}$ for $\delta=\delta(\sigma^3)
=E^{-5}[\iota_n,\eta^*]$. By the fact that\break 
$H(E^{-1}[\iota_{n-2},\eta\bar{\mu}])=4\bar{\zeta}$ (\fullref{T10}(2)) and $[\iota_{n-2},4\nu^*]\in E^6\pi^{n-8}_{2n+7}$ (\fullref{nnu*}(2)), 
$E^2\delta\in P\pi^{2n-7}_{2n+12}=0$, $E\delta\in\{[\iota_{n-4},4\bar{\kappa}]\}=0$ and ${\mathcal CDR}[\delta\in P\pi^{2n-9}_{2n+12}]$. 

Let $n\equiv 1\ (8)$. \fullref{nT22}(1) and 
${\mathcal ASM}[\nu^*]$ imply
$E^2\delta\in\{[\iota_{n-1},4\bar{\zeta}], 
[\iota_{n-1},\bar{\sigma}]\}\subset E^4\pi^{n-5}_{2n+12}$ (\fullref{nTa1}(1)), where $\delta=\delta(\sigma^3)=E^{-3}[\iota_n,\nu^*]$.  
Hence, $E\delta\in P\pi^{2n-3}_{2n+17}=0$ and 
${\mathcal CDR}[\delta\in P\pi^{2n-7}_{2n+14}]$. 

Let $n\equiv 0\ (8)$. 
\fullref{nnu*}(2) and ${\mathcal ASM}[4\nu^*]$ imply $E^5\delta\in P\pi^{2n-1}_{2n+18}=0$ and  
$E^4\delta\in\{[\iota_{n-2},4\bar{\kappa}]\}=0$ \eqref{nk6} for $\delta=\delta(\eta\eta^*\sigma)=E^{-6}[\iota_n,4\nu^*]$. Therefore, by \eqref{sgm35}, 
$E^3\delta\in\{E(\tau_{n-4}\eta\bar{\kappa})\}$ and  $E^2\delta\in\{[\iota_{n-4},\eta^2\bar{\kappa}], [\iota_{n-4},\nu\bar{\sigma}]\}$. 
By the relation $H(E^{-1}[\iota_{n-4},\eta^2\bar{\kappa}])=4\nu\bar{\kappa}$ (\fullref{T10}(2)) and \eqref{nbs4}, $E\delta\in\{[\iota_{n-5},\nu\bar{\kappa}], 
[\iota_{n-5},\bar{\rho}]\}$ $\subset E^3\pi^{n-8}_{2n+9}$, 
$\delta\in P\pi^{2n-11}_{2n+13}$ and ${\mathcal CDR}[\eta\eta^*\sigma=0]$.
\end{proof}

By Propositions \ref{T10}(4) and \ref{nu**}, 
$$
[\iota_{8n+1},\eta\eta^*]\neq 0.
$$
We show: 
\begin{prop}\label{et***}
$\sharp[\iota_{16n+14},\eta^*]=\sharp[\iota_{16n+13},\eta\eta^*]=2$.
\end{prop}
\begin{proof}
We use \fullref{Nnnu*}(4). Let $n\equiv 14\ (16)$. By
\fullref{nT22}(4), $[\iota_{n-1},\eta^2\rho]$ desuspends seven
dimensions. So, by the relation
$[\iota_{n-1},\bar{\mu}]=E(\tau_{n-2}\bar{\mu})$, \eqref{nukap5} and
${\mathcal ASM}[\eta^*]$,
$E^5\delta\in\{[\iota_{n-2},4\nu^*],[\iota_{n-2},\eta\bar{\mu}]\}$ for
$\delta=\delta(\eta^*\sigma)=E^{-7}[\iota_n,\eta^*]$. By the relation
$H(E^{-1}[\iota_{n-2},\eta\bar{\mu}])=4\bar{\zeta}$ and
\fullref{nnu*}(2),
$E^4\delta\in\{[\iota_{n-3},\bar{\zeta}],[\iota_{n-3},
\bar{\sigma}]\}$.  By the relation $\nu\bar{\zeta}=0$ and
\fullref{nnu*}(1), $E^3\delta\ (\bmod\
E^2(\bar{\tau}_{n-6}\bar{\sigma}))\in\{[\iota_{n-4},4\bar{\kappa}]\}$.
By \eqref{sgm39}, $E^2\delta\ (\bmod\
E(\bar{\tau}_{n-6}\bar{\sigma}),E^2\delta_1)\in\{E(\tau_{n-6}\eta\bar{\kappa})\}$,
where
$\delta_1=\delta(4\nu\bar{\kappa})=[\iota_{n-4},4\bar{\kappa}]$. This
induces $E\delta\ (\bmod\ E\delta_1)\in
P\pi^{2n-11}_{2n+11}=\{E\delta_2, [\iota_{n-6}, \nu\bar{\sigma}]\}$,
where $E\delta_2=[\iota_{n-6},\eta^2\bar{\kappa}]$,
$H\delta_2=4\nu\bar{\kappa}$ and $[\iota_{n-6},\nu\bar{\sigma}]\subset
E^2\pi^{n-8}_{2n+7}$ (\fullref{T10}(1)). Hence, $\delta\ (\bmod\
\delta_1,\delta_2)\in P\pi^{2n-13}_{2n+10}$ and ${\mathcal
CDR}[\eta^*\sigma\in 2\pi^s_{23}]$.

Next, let $n\equiv 13\ (16)$. ${\mathcal ASM}[\eta\eta^*]$ implies $E^6\delta\in\{[\iota_{n-1},4\nu^*],[\iota_{n-1},\eta\bar{\mu}]\}$ for $\delta=\delta(\eta\eta^*\sigma)=E^{-7}[\iota_n,\eta\eta^*]$. By the relation
$H(E^{-1}[\iota_{n-1},\eta\bar{\mu}])=\eta^2\bar{\mu}$ and \fullref{nnu*}(2),  
$E^5\delta\in\{[\iota_{n-2},\bar{\zeta}],[\iota_{n-2},\bar{\sigma}]\}$. By Lemmas \ref{nnu*}(1) and \ref{Nnnu*}(3), 
$E^4\delta\ (\bmod\ E^2(\bar{\tau}_{n-5}\bar{\sigma}))\in\{[\iota_{n-3},4\bar{\kappa}]\}\subset E^6\pi^{n-9}_{2n+7}$ and $E^3\delta\ (\bmod\ E(\bar{\tau}_{n-5}\bar{\sigma}))\in\{[\iota_{n-4},\eta\bar{\kappa}],[\iota_{n-4},\sigma^3]\}$. From the fact that $[\iota_{n-4},\eta\bar{\kappa}]=E(\tau_{n-5}\eta\bar{\kappa})$ and \eqref{sgm39},  $E^2\delta\in\{[\iota_{n-5},\eta^2\bar{\kappa}],
[\iota_{n-5},\nu\bar{\sigma}]\}$. 
Since $H(E^{-1}[\iota_{n-5},\eta^2\bar{\kappa}])=
4\nu\bar{\kappa}$ and
$H(E^{-3}[\iota_{n-5},\nu\bar{\sigma}])
=\nu^2\bar{\sigma}=0$ (\fullref{nTa1}(1),\cite{Od2}), 
$E\delta\in P\pi^{2n-11}_{2n+12}\subset E^7\pi^{n-13}_{2n+3}$ \eqref{des7}, $\delta\in\ P\pi^{2n-13}_{2n+11}$ and ${\mathcal CDR}[\eta\eta^*\sigma=0]$. 
\end{proof} 

We show the following:
\begin{prop}\label{bsigma1}
$H(E^{-3}[\iota_{8n+k},\bar{\sigma}]_{\neq 0})=\nu\bar{\sigma}$ for $k=0,1,2$. 
\end{prop}
\begin{proof}
Let $n\equiv 0\ (8)$. By Lemmas \ref{nTa1}(1), \ref{nT11} and \ref{nT22}, there exists an element $\delta(k)\in\pi^{n+k-3}_{2n+2k+15}$ such that $[\iota_{n+k},\bar{\sigma}]=E^3\delta(k)$ and $H\delta(k)=\nu\bar{\sigma}$. For $k=0$, ${\mathcal ASM}[\bar{\sigma}]$ induces $E^2\delta(0)\in P\pi^{2n-1}_{2n+19}=0$, $E\delta(0)\in P\pi^{2n-3}_{2n+18}\subset E^2\pi^{n-4}_{2n+14}$ (\fullref{T11}(1)) and ${\mathcal CDR}[\delta(0)\in P\pi^{2n-5}_{2n+17}]$. 
 By the parallel argument to \fullref{nu**} for $\nu^*$, the assertion follows for $k=1$. For $k=2$, ${\mathcal ASM}[\bar{\sigma}]$ induces  $E^2\delta(2)\in\{E(\tau_n\bar{\kappa})\}$ and  
$E\delta(2)\in\{[\iota_n,\eta\bar{\kappa}],[\iota_n,\sigma^3]\}$. Since $[\iota_n,\sigma^3]\subset E^7\pi^{n-7}_{2n+13}$  (\fullref{nTa1}(3)), we obtain $\delta(2)\ (\bmod\ \beta)=0$ and  
${\mathcal CDR}[\nu\bar{\sigma}\ (\bmod\ \eta^2\bar{\kappa})=0]$, where $\beta=\delta(\eta^2\bar{\kappa})=E^{-1}[\iota_n,\eta\bar{\kappa}]$ (\fullref{T10}(1)). 
\end{proof}

We show the following:
\begin{prop}\label{etbkap}
$H(E^{-5}[\iota_{8n+2},\eta\bar{\kappa}]_{\neq 0})=\nu^2\bar{\kappa}$ and $H(E^{-6}[\iota_{8n+1},
\eta^2\bar{\kappa}]_{\neq 0})
=\varepsilon\bar{\kappa}$.
\end{prop}
\begin{proof}
Let $n\equiv 2\ (8)$. By \fullref{No}(1) and \eqref{bb}, 
$H(E^{-5}[\iota_n,\eta\bar{\kappa}])
=\nu^2\bar{\kappa}.$ We set $\delta=\delta(\nu^2)=E^{-5}[\iota_n,\eta]$. 
${\mathcal ASM}[\eta\bar{\kappa}]$ induces  
$E^4(\delta\bar{\kappa})\in P\pi^{2n-1}_{2n+21}\subset 
E^5\pi^{n-6}_{2n+14}$ (Lemmas \ref{nT22}(3),\ref{No}(2)) and  
$E^3(\delta\bar{\kappa})\in\{[\iota_{n-2},4\nu\bar{\kappa}], [\iota_{n-2},
8\bar{\rho}]$, 
$[\iota_{n-2},\eta^*\sigma]\}$. By \fullref{nTa1}(1), the first two 
Whitehead products desuspend four dimensions, respectively. Hence, by the relation $H(E^{-1}[\iota_{n-2},\eta^*\sigma])=\eta\eta^*\sigma$, we obtain 
$E^2(\delta\bar{\kappa})=0$,  
$E(\delta\bar{\kappa})\in P\pi^{2n-7}_{2n+18}\subset E^2\pi^{n-6}_{2n+14}$ (\fullref{T11}(1)), $\delta\bar{\kappa}\in P\pi^{2n-9}_{2n+17}$ and  ${\mathcal CDR}[\nu^2\bar{\kappa}=0]$. 

Next, let $n\equiv 1\ (8)$.
By \fullref{No}(2),  
$H(E^{-6}[\iota_n,\eta^2\bar{\kappa}])
=\varepsilon\bar{\kappa}$. 
${\mathcal ASM}[\eta^2\bar{\kappa}]$ implies 
$E^5(\delta\bar{\kappa})\in\{[\iota_{n-1},4\nu\bar{\kappa}],[\iota_{n-1},8\bar{\rho}],[\iota_{n-1},\eta^*\sigma]\}$ for $\delta=\delta(\varepsilon)=E^{-6}[\iota_n,\eta^2]$. 
By \fullref{nTa1}(2), $[\iota_{n-1},8\bar{\rho}]$ desuspends eight dimensions. By \fullref{nTa11}(3), 
$[\iota_{n-1},4\nu\bar{\kappa}]=[\iota_{n-1},\nu^3]\kappa$ desuspends six dimensions. So, by the relation\break  
$H(E^{-1}[\iota_{n-1},\eta^*\sigma])=\eta\eta^*\sigma$,
we have $E^4(\delta\bar{\kappa})\in
P\pi^{2n-3}_{2n-27}=0$,  
$E^3(\delta\bar{\kappa})\in\{[\iota_{n-3},\mu_{3,\ast}]\}$
$\subset E^6\pi^{n-9}_{2n+12}$ (\fullref{6bamu}), $E^2(\delta\bar{\kappa})\in\{[\iota_{n-4},\eta\mu_{3,\ast}]\}\subset E^4\pi^{n-8}_{2n+13}$ (\fullref{nT22}(1)),  
$E(\delta\bar{\kappa})\in\{[\iota_{n-5},4\zeta_{3,\ast}]\}\subset E^3\pi^{n-8}_{2n+7}$ 
(\fullref{T11}(3)), $\delta\bar{\kappa}\in P\pi^{2n-11}_{2n+17}$ and hence, ${\mathcal CDR}[\varepsilon\bar{\kappa}=0]$. 
\end{proof}

According to Mahowald \cite{M3}, the following seems to be true.

\begin{con} \label{Mah}
$\langle\nu,\eta,\bar{\sigma}\rangle
=\langle\bar{\nu},\sigma,\bar{\nu}\rangle
=\eta\eta^*\sigma$.
\end{con}

By use of the Jacobi identity for Toda brackets, 
\fullref{Mah} and the relations $\langle\eta,\nu,\eta\rangle=\nu^2,\sigma\bar{\sigma}=0$ \cite{Od2}, we obtain 
$$
\langle 2\iota,\nu^2,\bar{\sigma}\rangle=\langle 2\iota,\eta,\eta\eta^*\sigma\rangle=\nu^2\bar{\kappa}.
$$
By this fact, we can show 
$$
[\iota_{8n},\nu\bar{\sigma}]\neq 0.
$$
\begin{proof} 
Let $n\equiv 0\ (8)$. In \fullref{nTkl}[n-5;5,3], $\P^{n-1}_{n-5}=E^{n-8}\P^7_3$ and $\gamma_{7,3}=2s_4 + i^7_3\tilde{\eta}''$, where $s_4=p^7_3s_3$ \eqref{gam7}. By \fullref{coexnu2}(8),
$$
\tilde{\eta}''\circ\nu\bar{\sigma}\in i^{4,6}_3\circ\langle i\bar{\eta},\tilde{\eta},\nu\rangle\circ\bar{\sigma}
=i\circ\langle 2\iota,\nu^2,\bar{\sigma}\rangle=
i\nu^2\bar{\kappa}.
$$
This shows 
$$
H(E^{-4}[\iota_n,\nu\bar{\sigma}])=\nu^2\bar{\kappa}.
$$
For $\delta=\delta(\nu^2\bar{\kappa})=E^{-4}[\iota_n,\nu\bar{\sigma}]$, 
${\mathcal ASM}[\nu\bar{\sigma}]$ implies 
$E^3\delta=0$ and $E^2\delta\in P\pi^{2n-3}_{2n+21}$\\
$\subset E^3\pi^{n-5}_{2n+16}$ (\fullref{T11}(1)), $E\delta\in\{[\iota_{n-3},\eta^2\bar{\rho}],[\iota_{n-3},\mu_{3,\ast}]\}$, \\ 
$\delta\ (\bmod\ \tau_{n-4}\eta^2\bar{\rho},
\tau_{n-4}\mu_{3,\ast})\in P\pi^{2n-7}_{2n+19}$ and hence, 
${\mathcal CDR}[\nu^2\bar{\kappa}\ (\bmod\ \eta\mu_{3,\ast})=0]$. 
\end{proof}

Finally, by \fullref{T11}(1) and \fullref{nT22}(1), we note the following. 

\begin{rem}
 $H(E^{-2}[\iota_{8n+2},4\bar{\kappa}])
=\varepsilon\kappa
=\eta^2\bar{\kappa}$ and 
$H(E^{-3}[\iota_{8n+1},\nu\bar{\kappa}])
=\nu^2\bar{\kappa}$.
\end{rem}

\bibliographystyle{gtart}
\bibliography{link}

\end{document}